\documentclass[11pt,a4paper]{article}
\usepackage[left=3cm, right=3cm, top=3cm, bottom=3cm]{geometry}
\usepackage{graphicx}
\usepackage[T1]{fontenc}
\usepackage[utf8]{inputenc}
\usepackage[overload]{empheq}
\usepackage{setspace}
\usepackage{amsmath}
\usepackage{amssymb}
\usepackage{amsfonts}
\usepackage{relsize}
\usepackage{bm}
\usepackage[dvipsnames]{xcolor}
\usepackage{enumitem}
\usepackage{url}
\usepackage{pdfpages}
\usepackage{palatino}
\usepackage{wrapfig}

\makeatletter
\def\patchsect#1\let\@svsec\@empty{#1\def\@svsec{\leavevmode\kern1sp\relax}}
\let\old@sect\@sect
\def\@sect{\expandafter\patchsect\old@sect}
\makeatother

\usepackage{amsthm}
\theoremstyle{definition}

\theoremstyle{plain}
\newtheorem{defi}{Definition}[section]
\newtheorem{lemma}[defi]{Lemma}
\newtheorem{theorem}[defi]{Theorem}
\newtheorem{prop}[defi]{Proposition}

\newtheorem{remark}[defi]{Remark}

\usepackage{titlesec}
\titleformat{\section}
  {\centering\normalfont\large\bfseries}{\thesection.}{1em}{}
\titleformat{\subsection}
  {\normalfont\large\bfseries}{\thesubsection.}{1em}{}
\titleformat{\subsubsection}
  {\normalfont\normalsize\itshape\color{MidnightBlue}}{\thesubsubsection.}{1em}{}

\usepackage{hyperref}
 \hypersetup{
    colorlinks,%
    citecolor=red,%
    filecolor=black,%
    linkcolor=black,%
    urlcolor=black
}

\usepackage{float}

\usepackage{subcaption}
\usepackage{icomma}
\usepackage{appendix}
\usepackage[backend=bibtex, backref=true, style=numeric,maxbibnames=99, url=false, isbn=false]{biblatex}

\bibliography{references}
\DeclareMathOperator*{\Id}{Id}
\newcommand{\jap}[1]{\langle #1 \rangle}

\title{Weak-strong uniqueness for the Landau equation by a relative entropy method}
\date{\today}
\author{Côme Tabary\footnote{DMA, Ecole Normale Supérieure, Paris}}
\begin{document}

\maketitle

\begin{abstract}
\noindent We derive a weak-strong uniqueness and stability principle for the Landau equation in the soft potentials case (including Coulomb interactions). The distance between two solutions is measured by their relative entropy, which to our knowledge was never used before in stability estimates. The logarithm of the strong solution is required to have polynomial growth while the weak solution can be any H-solution with sufficiently many moments at initial time. Since we require a substantial amount of regularity on the strong solution, we also provide an example of sufficient conditions on the initial data that ensure this regularity in the Coulomb (and very soft potentials) case.
\end{abstract}

\section{Introduction}

\subsection{Background and main results}

In this work, we consider the space-homogeneous Landau equation in three-dimensional velocity space:
\begin{equation}
\label{eq:landau}
    \partial_t f_t(v) = \nabla_v \cdot \int_{\mathbb{R}^3} a (v - v') (\nabla_v - \nabla_{v'}) f(v)f(v') dv',
\end{equation}
with the matrix-valued function $a(z)=\vert z \vert^{\gamma+2} \Pi(z)$, where $$\Pi(z) = \left(I - \frac{z}{\vert z \vert} \otimes \frac{z}{\vert z \vert}\right)$$
is the projection on $z^\perp$. The unknown $f_t:\mathbb{R}^3\mapsto \mathbb{R}_+$ models the distribution of velocities in a plasma at time $t$. Obeying basic physical principles, the Landau equation preserves the total mass, momentum and energy of the plasma. Thanks to the mass conservation, up to rescaling the time variable, we can and will consider distributions of velocities with normalized unit mass, \textit{i.e.} probability densities:
$$\forall t\geq0,\ \int_{\mathbb{R}^3}f_t =1.$$
The (mathematical) entropy
$$H(f)=\int_{\mathbb{R}^3} f \ln f$$ decreases along trajectories. Formally,
$$\frac{d}{dt}H(f_t) = -\frac{1}{2}\iint_{\mathbb{R}^6}a (v - v') :\left[(\nabla - \nabla')\ln(f_tf_t') \right]^{\otimes 2}f_tf_t'dvdv' =:-D(f_t)$$
which is non-positive because $a$ is non-negative.
The Landau equation has been widely used in plasma physics as it can be derived from the Boltzmann equation of gas dynamics in the so-called \textit{grazing collisions limit}. Moreover, in the case of Coulomb interactions, corresponding to $\gamma=-3$, the Landau equation is the reasonable alternative to the Boltzmann equation which becomes meaningless \cite{LifsicPitaevskij2008}. The most physically relevant case hence being $\gamma=-3$, we concentrate on \textit{very soft potentials} $\gamma\in[-3,-2]$ for applications. However, our principal result holds for the full range $\gamma\in[-3,0)$ of soft potentials. 

The main result of this work is a weak-strong uniqueness and stability estimate, under conditional logarithmic bounds on the strong solution. The weak formulation we use in this work is Villani's notion of H-solutions \cite{Villani1998} (see Definition \ref{def:wsol} below for a precise definition).
\begin{theorem}
\label{thm:rentropy}
    Fix $\gamma\in [-3,0)$ and a final time $T>0$. Let $g=(f_t)_{t\in[0,T]}$ be a H-solution of the Landau equation on $[0,T]$ with initial data $f_0$. Let $g=(g_t)_{t\in[0,T]}$ be a classical solution to the Landau equation for $C^2$ initial data $g_0>0$ (with unit mass, finite energy and finite entropy). Further assume that $g$ satisfies the following:
    \begin{itemize}
        \item (Qualitative smoothness) $g$ is $C^1$ in time and $C^2$ in velocity.
        \item (Logarithmic tail control) It admits the following bounds: for some $\kappa,\nu \geq 1$,
    $$\left\Vert \jap{v}^{-\kappa} \nabla \ln g_s \right\Vert_{L^2([0,T],L^\infty(\mathbb{R}^3))}<+\infty,$$
    $$\left\Vert \jap{v}^{-\nu} \partial_s \ln g_s \right\Vert_{L^1([0,T],L^\infty(\mathbb{R}^3))}<+\infty.$$
     If $\gamma<-2$, for some $\zeta\geq 0$,
    $$\left\Vert \jap{v}^{-\zeta} \nabla^2 \ln g_s \right\Vert_{L^2([0,T],L^\infty(\mathbb{R}^3))}<+\infty,$$
        \item  (Moments) Assume that for $\rho=\max(2\kappa, 2\kappa+2\gamma +4,2\zeta,\nu)$, the $(\rho-\gamma)$-th moment of $g_0$ is bounded:
    $$ \int g_0(v)\jap{v}^{\rho-\gamma}dv= M_{g_0} <\infty.$$
    \end{itemize}
    Further assume that the $(\rho-\gamma)$-th moment of $f_0$ is also bounded:
    $$\int f_0(v)\jap{v}^{\rho-\gamma}dv = M_{f_0} <\infty.$$
    Then the relative entropy satisfies for all times $t\in[0,T]$:
    \begin{equation}
    \label{eq:thm_rentropy}
    H(f_t \vert g_t) \leq  H(f_0 \vert g_0) \exp\left(\mathcal{C}\int_0^t \left(\Vert \jap{v}^{-\kappa} \nabla \ln g_s \Vert_{L^\infty(\mathbb{R}^3)} ^2+\Vert \jap{v}^{-\zeta} \nabla^2 \ln g_s \Vert_{L^\infty(\mathbb{R}^3)} ^2\right)ds\right)  
    \end{equation}
    with $\mathcal{C}>0$ a constant depending on $\gamma,T,\rho,M_{g_0},H(g_0),M_{f_0},H(f_0)$.
\end{theorem}

The above conditions on $g$ are quite strong and not very usual in the literature (although logarithmic bounds for the Landau equation with $\gamma=0$ were also needed by the authors of \cite{CarrilloFengGuo2024} to close their relative entropy estimate).
Note that the relative entropy estimate \eqref{eq:thm_rentropy} does not depend on $\Vert \jap{v}^{-\nu}\partial_s \ln g_s \Vert_{L^1([0,T], L^\infty)}$, which could be an indication that this hypothesis can be relaxed. It is indeed quite (too) strong, because integrating it yields $L^\infty$ in time control on $\ln g$ which restricts the initial data. We only use it to pass to the limit in the weak formulation of H-solutions, which might be done by a better approximation argument.

In the very soft potential case $\gamma \in [-3,-2]$, we present below sufficient conditions on the initial data for a strong solution satisfying the conditions of Thereom \ref{thm:rentropy} to exist. Its statement involves Hölder continuity and associated seminorms, their notation is recalled in Section \ref{ssec:definitions}.

\begin{theorem}
\label{thm:gbounds}
Let $\gamma \in [-3,-2]$. Let $g_0\in L^\infty(\mathbb{R}^3)$ be a initial data with unit mass, finite energy and entropy. Further suppose that
\begin{itemize}
    \item (Tail control) $g_0$ is bounded above and below by a multiple of the same Maxwellian:
    \begin{equation}
    \label{eq:maxwellcond}
         ke^{-\mu\vert v \vert^2}\leq g_0(v) \leq Ke^{-\mu\vert v \vert^2},
    \end{equation}
    for some $\mu,k,K>0$.
    \item (Smoothness) The function $\ln g_0$ is locally $C^{2+\delta}$ for some $\delta \in (0,1)$, and there exists $\ell\in \mathbb{R}$ such that, for any $R\geq 1$,
    \begin{equation}
    \label{eq:smoothcond}
         \sup_{i,j}[\partial_{ij} \ln g_0]_{\delta,B(0,R)} \leq C_g R^\ell
    \end{equation}
    for some constant $C_g>0$.
    
\end{itemize}
Under these conditions, there exists a (unique) classical solution $g=(g_t)_{t\in[0,+\infty)}$ of the Landau equation \eqref{eq:landau}.

Let $f_0 \in L^1(\mathbb{R}^3)$ be a non-negative initial data with unit mass, finite energy and entropy, and $f=(f_t)_{t\in[0,+\infty)}$ an associated weak solution to the Landau equation. There exists $k=k(\gamma,\ell,\delta)>0$ such that, if $\jap{v}^kf_0 \in L^1$, then the following relative entropy estimate holds: for all $t\in[0,+\infty)$,
$$H(f_t\vert g_t) \leq C(t) H(f_0 \vert g_0),$$
where $C(t)$ depends only on $\gamma, f_0, g_0$ and $t$.
\end{theorem}
\begin{remark}
The condition \eqref{eq:smoothcond} directly implies that $g_0$ has finite energy (in fact, infinitely many moments) and finite entropy.
\end{remark}

\begin{remark}
    The exact expression of the moment is $$k(\gamma,\ell,\delta)=8\left(\max\left(2+\frac{1}{2+\delta}(\ell-2), 1\right)+\frac{2-\gamma}{2\delta}\right)-\gamma.$$ If $\gamma=-3$, $\ell \leq -1$ and $\delta$ is arbitrary close to $1$ (this holds for $g_0$ of class $C^3$ with Maxwellian tails), the initial condition of the weak solution $f_0$ must have a moment of order strictly greater than $31$.
\end{remark}

We believe the above result could be extended to the less physically relevant case $\gamma>-2$ through other means but did not investigate this further.
We prove Theorem~\ref{thm:gbounds} by showing that the strong solution satisfies all the necessary bounds for Theorem \ref{thm:rentropy} to apply. More precisely, we actually show in the proof that the bounds in velocity hold in $L^\infty$ in time, provided they hold at $t=0$ (and with additional regularity on the initial data). It is reassuring that the regularity we ask in Thereom \ref{thm:rentropy} is propagated by the Landau equation (otherwise Thereom \ref{thm:rentropy} would be quite meaningless in practice), however it would be preferable to \textit{generate} these bounds for positive times to relax the assumptions on the initial data. This question is beyond the scope of this first result and we postpone it to a future work.

\subsection{Review of the literature and discussion}

The literature on the Landau equation is expansive. Far from being exhaustive, we present below the main results that motivated this work.
\\

\noindent\textbf{Weak formulations and uniqueness.} The major difficulty in the study of the very soft potentials is the singularity of $a$ at $z=0$. The straightforward weak formulation of the right hand side in \eqref{eq:landau}, obtained by testing against a smooth function and integrating by parts, still requires $f$ to lie in some $L^p$ space (in velocity) to make sense of the integration against the singularity (typically, $f\in L^1([0,T],L^{p})$ for some $p>\frac{3}{\gamma+5}$). Some even stronger integrability of $f$ is also required for the uniqueness results of Fournier and Guerin \cite{Fournier2008, GuerinFournier2008}: there is at most one weak solution in $L^1([0,T], L^{3/({\gamma+3})})$.

Villani bypassed the integrability condition for existence of solutions by introducing the notion of H-solutions \cite{Villani1998}. Using the production of entropy $D(f_t)$, one can make sense of
the integral in \eqref{eq:landau} in a \textit{a priori} even weaker way. The production can indeed be formally rewritten
$$D(f_t)=2\iint_{\mathbb{R}^6}\vert v-v'\vert^{\gamma+2} \left\vert \Pi (v - v') (\nabla - \nabla')\sqrt{f_tf_t'} \right\vert^{2}dvdv', $$
becoming a weighted $L^2$ norm of the (weak) gradient of $\sqrt{ff'}$ in a particular direction, which is used to give sense to a weak formulation (see Definition \ref{def:wsol} and the following remark for details). Relying on $D(f_t)$ is more natural because the Landau equation directly provides an a priori estimate for the entropy production, rather than for an $L^p$ norm. However, Desvillettes later showed that the production of entropy controls a weighted $L^3$ norm of $f_t$, granting enough integrability to H-solutions to ensure they are weak solutions in the more classical sense \cite{Desvillettes2014a}. The weight was recently improved by Ji in \cite{Ji2024b}, where it is also shown that the entropy dissipation cannot control a $L^p$ norm with $p>3$. In the $\gamma\in[-3,-2)$ regime, this means one cannot reach the integrability required to apply the uniqueness theorem from \cite{Fournier2008, GuerinFournier2008} using the entropy production alone.

Another result on uniqueness is the more recent work by Chern and Gualdani~\cite{ChernGualdani2022}, where uniqueness is shown for weak solutions in $L^\infty([0,T], L^p)$ for some $p>\frac{3}{2}$, and satisfying the entropy production estimate (and some qualitative regularity which allows to use truncations of the solution as test functions). This was extended to more general $L^q([0,T], L^p)$ conditions through a Prodi-Serrin-type criterion \cite{AlonsoBaglandDesvillettes2024} (including the case $q=1$, $p=+\infty$ from \cite{Fournier2010}). We emphasize that uniqueness is only one of the byproducts and that the satisfaction of the Prodi-Serrin criterion yields much more regularity on the solutions, which makes it an important threshold between very weak solutions and more regular ones. However, as of current results, these integrability conditions cannot \textit{a priori} stem from the Landau equation itself.

Improving any of the aforementioned uniqueness results to work only with the regularity obtained from the entropy production alone seems quite a difficult task. Instead, we focus on weak-strong uniqueness, since it can cover the case of even weaker solutions because nearly all regularity assumptions are transferred to the strong one. Taking $f_0=g_0$ in Theorem \ref{thm:gbounds} for instance, we have a uniqueness result in the class of all H-solutions starting from $g_0$. In a way, we trade the integrability conditions from Fournier and Guerin \cite{Fournier2010,GuerinFournier2008} with the entropy production estimate of H-solutions (see \eqref{eq:entropyprod} below), which is both more physical and mathematically weaker, since it only implies $f\in L^1([0,T],L^3)$). However, the regularity of the initial data we demand is still too strong for our uniqueness result to be entirely comparable to \cite{AlonsoBaglandDesvillettes2024,ChernGualdani2022,Fournier2010,GuerinFournier2008}. The fact that whenever there exists a regular enough solution, any weak solution coincides with it is to be expected from any reasonable notion of weak solution, and it could be argued that, on physical grounds, such weak-strong uniqueness might in fact be more important than uniqueness of weak solutions in very irregular regimes where no strong solution exists at all.
\\

\noindent \textbf{Strong solutions in the class of weak solutions.} Weak-strong uniqueness is all the more interesting in view of the recent major advancements in the existence of smooth solutions, which now covers nearly any physically reasonable initial data. The breakthrough by Guillen and Silvestre \cite{GuillenSilvestre2023} on the monotonicity of the Fisher information for \eqref{eq:landau} shows existence of smooth classical solutions for smooth initial data with Maxwellian tails. This was then improved to initial data that is merely in $L^1$ (with a reasonable amount of moments) \cite{Ji2024, DesvillettesGoldingGualdani2024}.

These solutions are unique \textit{in the class of smooth solutions}, however nothing prevents the existence of a weaker solution, even for very well-behaved initial data. One could argue that once a smooth solution is shown to exist, there is no reason to enlarge the class of solutions to include weaker ones anymore. However, the question of wether uniqueness holds in the class of all weak solutions is still relevant for several reasons: First, because for rough initial data, even though the solutions are smooth for $t>0$, the initial data is still achieved in a weak sense so there could be a bifurcation at initial time, between two equally regular solutions (this is discussed in \cite{Ji2024}). Second, because the notion of weak solution is handier and more stable when taking limits or using approximation arguments: by taking the grazing collision limit of the Boltzmann equation or the limit $N\rightarrow\infty$ of some $N$-particle system, one usually recovers a weak solution. In this setting, one cannot guarantee that the weak solution is obtained by a suitable well-behaved approximation scheme and one has to work with what little regularity goes through in the (often weak) limit. A weak strong principle guarantees that the limit is in fact the strong solution without requiring a stronger form of convergence. Finally, we believe it is an interesting question in itself to figure the minimum regularity to be imposed to eliminate possible pathological solutions to the Landau equation: is the sole control of the physical quantities of mass, energy, entropy and entropy production enough? This work provides a partial positive answer. 
\\

\noindent \textbf{Relative entropy methods.}
Uniqueness and stability estimates in the literature are typically formulated in the Wasserstein distance (the aforementioned \cite{Fournier2008, GuerinFournier2008}, but also \cite{FournierHauray2015,NicolasFournierArnaudGuillin2017}) which corresponds to weak convergence of probability measures and additional convergence of some moments. We rather work with the relative entropy, which is mathematically better suited when one of the solutions is very regular, as it allows for explicit computations, but also because of physical motivations.

Indeed, in the pioneering works by Dafermos \cite{Dafermos1979} and Diperna \cite{Diperna1979}, a close connection between stability and the second principle of thermodynamics was exposed, in the setting of hyperbolic systems of conservation laws. It is natural to expect this connection to also exist at the level of the Landau equation, where the second principle manifests itself through the monotonicity of the entropy $H$. The relative entropy in our work also plays a similar role as the relative energy ($L^2$ norm) does in fluid dynamics (see \cite{Wiedemann2017} for a review, and Theorem 2.1 therein for a typical example of weak-strong uniqueness). In all weak-strong principles, the weak solution is asked to satisfy an \textit{admissibility criterion} which takes the form of an entropy production (for hyperbolic systems) or energy dissipation (for fluids) \textit{inequality}. For the Landau equation, this criterion is exactly the entropy production estimate embedded in the definition of H-solutions (see \eqref{eq:entropyprod} below).

To our knowledge, relative entropy/energy methods were rarely applied to the homogeneous Landau equation. We can however mention the recent work \cite{CarrilloFengGuo2024} in which the relative entropy between the solution to the Landau equation (for $\gamma=0$) and the distribution of an approximating particle system is used to show propagation of chaos.

The relative entropy provides stability results in strong topology since it controls the $L^1$ norm. To our knowledge, this is the first result of this kind in the literature. It is also interesting to find proofs of uniqueness that do not rely on a probabilistic interpretation of the Landau equation. 
\\

\textbf{Regularity estimates for the Landau equation.} Numerous works have studied the regularity of solutions to the Landau equation in the range of soft potentials. To apply Theorem \ref{thm:rentropy}, we are interested in estimating the $C^2$ regularity of classical solutions. The approach to obtain such higher regularity typically consists in writing the Landau equation \eqref{eq:landau} as a parabolic equation with coefficients depending on the solution, derive ellipticity estimates for the coefficients and apply Schauder estimates: see for instance \cite{HendersonSnelson2019, HendersonSnelsonTarfulea2019, HendersonWang2024, ImbertMouhot2021} (in the more general inhomogeneous setting). Without regularity assumptions on the initial data, these estimates of course degenerate as $t\rightarrow 0$ as powers of $\frac{1}{t}$. These powers are typically not square-integrable in time, so we do not obtain the $L^2$ bounds in time that are needed to apply Theorem \ref{thm:rentropy}. A similar discussion on the issue of time-integrability of estimates was done in \cite[Section 1.4.2]{HendersonSnelsonTarfulea2019}. It thus seems necessary to suppose \textit{some} regularity at time $t=0$. We deal with this issue by boldly assuming that at time $t=0$ we basically already have the regularity needed and propagate it for $t>0$: we show in the proof that the logarithmic tail control in Theorem~\ref{thm:rentropy} holds in $L^\infty$ in time. This is of course not optimal and there should exist a less restrictive assumption on the initial data that provides estimates that blow up at $t=0$ but that remain square-integrable. Since going from estimates on the derivatives of $g$ to ones on the derivatives of $\ln g$ already requires some technical work, we did not investigate further those potential improvements on the smoothness hypothesis of Theorem \ref{thm:gbounds}.

\subsection{Definitions and notation}
\label{ssec:definitions}
We gather here the notation used in this work. We measure the closeness of two solutions through their relative entropy: for two probability densities $f$ and $g$, we define
$$H(f\vert g) :=  \int_{\mathbb{R}^3} f\ln\left(\frac{f}{g}\right)=\int_{\mathbb{R}^3} \left(\frac{f}{g}\ln\left(\frac{f}{g}\right)-\frac{f}{g}+1\right)g.$$
Note that the right hand side always has a value in $[0,+\infty]$ as it is the integral of a non-negative function.
The relative entropy is known to control the (square root of the) $L^1$ norm by Pinsker's inequality.

We recall the definition of H-solutions (as in \cite{Villani1998}):
\begin{defi}
\label{def:wsol}
Given $\gamma\in[-3,0)$ and initial data $f_0\in L^1(\mathbb{R}^3)$ with unit mass, finite energy (second moment) and finite entropy, a \emph{H-solution} to the Landau equation on $[0,T]$ is a non-negative function $f\in C([0,T],\mathcal{D}'(\mathbb{R}^3))\cap L^\infty([0,T],L^1(\mathbb{R}^3))$ satisfying for all $t\in[0,T]$ the conservation of mass, momentum and energy
\begin{align*}
\int f_t dv &= \int f_0 dv=1\\
\int f_t v dv &= \int f_0v dv\\
\int f_t \vert v \vert^2 dv &= \int f_0 \vert v \vert^2 dv,
\end{align*}
the entropy production estimate
\begin{equation}
\label{eq:entropyprod}
H(f_t) +2\int_0^t\iint_{\mathbb{R}^6}a (v - v') :\left[(\nabla - \nabla')\sqrt{f_sf_s'} \right]^{\otimes 2}dvdv'ds \leq H(f_0)
\end{equation}
as well as the following weak formulation of the equation: for any $\varphi\in C^2_0([0,T]\times \mathbb{R}^3)$, any $t\in[0,T]$:
    \begin{align*}
    \label{eq:testg}
     \int_{\mathbb{R}^3} \varphi_tf_tdv-\int_{\mathbb{R}^3} \varphi_0f_0dv &= \int_0^t \int_{\mathbb{R}^3} (\partial_s\varphi_s)f_s dv\ ds \\
    & -\int_0^t \iint_{\mathbb{R}^6} a(v-v') (\nabla\varphi_s - \nabla'\varphi'_s)\cdot\left((\nabla - \nabla')\sqrt{f_sf'_s}\right) \sqrt{f_s f_s'}\ dvdv'ds.
    \end{align*}
\end{defi}

\begin{remark}
\label{rem:weaksense}
The entropy production estimate implies that the quantity
$$ \int_0^T \iint_{\mathbb{R}^6} \vert v-v'\vert^{\gamma +2} \left\vert \Pi (v - v')(\nabla - \nabla')\sqrt{f_s f'_s} \right\vert^2 dvdv'$$
is finite, which allows one to make sense of the weak formulation: by Cauchy-Schwarz,
    \begin{align*}
&\int_0^t\iint_{\mathbb{R}^6}  \vert v-v'\vert^{\gamma +2}\Pi(v-v')(\nabla\varphi_s - \nabla'\varphi'_s)\cdot\left((\nabla - \nabla')\sqrt{f_sf'_s}\right)  \sqrt{f_sf'_s}dvdv'ds\\
&\leq \left(\int_0^t\iint_{\mathbb{R}^6}\vert v-v'\vert^{\gamma +2}  f_sf'_s \vert\nabla\varphi_s - \nabla'\varphi'_s\vert^2dvdv'ds\right)^\frac{1}{2}\\
&\ \; \; \; \cdot\left(\int_0^t\iint_{\mathbb{R}^6}\vert v-v'\vert^{\gamma +2}\left\vert\Pi(v-v')(\nabla - \nabla')\sqrt{f_sf'_s} \right\vert^2 dvdv'ds\right)^\frac{1}{2}.\\
\end{align*}
The second integral is finite by the entropy production estimate. If $\gamma+2<0$, we must deal with the singularity $v=v'$ in the first integral. In this case, the smoothness of $\nabla\varphi_s$ entails that
$\vert\nabla\varphi_s - \nabla'\varphi'_s\vert \leq C_\varphi \min(\vert v-v'\vert,1)$, \textit{i.e.} it is Lipschitz and bounded. In particular $\vert\nabla\varphi_s - \nabla'\varphi'_s\vert^2 \leq C_\varphi \vert v-v'\vert^{-\gamma-2}$, so the singularity is canceled. If $\gamma+2\geq0$, the bound $\vert\nabla\varphi_s - \nabla'\varphi'_s\vert^2\vert v-v'\vert^{\gamma +2} \leq C_\varphi (1+ \vert v \vert^2 + \vert v' \vert^2)$ holds and the energy conservation then ensures that the integral makes sense.
\end{remark}

We also make precise our notation for Hölder seminorms, which will be used in Section~\ref{sec:logarithmicbounds}:
for a function $\varphi$ defined on a domain $Q\subset \mathbb{R}_+\times\mathbb{R}^3$ and $\alpha\in(0,1)$, we define the seminorm
$$[\varphi]_{\alpha,Q} := \sup_{\substack{(t,v), (t,w) \in Q\\
v\neq w}} \frac{\vert \varphi(t,v)-\varphi(t,w)\vert}{\vert v-w \vert^\alpha},$$
which measures Hölder continuity \textit{in velocity only}, uniformly in time. If $Q\subset\mathbb{R}^3$ so that there is no time dependence, $[\varphi]_{\alpha,Q}$ is the usual Hölder seminorm.
For a differentiable or twice differentiable function, we also let
$$[\varphi]_{1+\alpha,Q}:=\sup_{i}[\partial_{i}\varphi]_{\alpha,Q}$$
$$[\varphi]_{2+\alpha,Q}:=\sup_{ij}[\partial_{ij}\varphi]_{\alpha,Q}$$
and say that $\varphi$ is $C^{k+\alpha}$ if $[\varphi]_{k+\alpha,Q}<+\infty$, for $k=0,1,2$. Interpolation results between these seminorms are recalled in Appendix~\ref{app:interpol}.

We use the japanese bracket $\jap{v}=\sqrt{1+\vert v \vert^2}$. A function $\varphi$ is said to admit a moment of order $k\in\mathbb{R}$ (or a $k$-th moment) if $(\jap{v}^k\varphi )\in L^1$.

Euclidean balls of radius $R>0$ centered at $v\in\mathbb{R}^3$ are denoted by $B(v,R)$.

The matrix dot product $A:B$ is $\sum_{ij} A_{ij} B_{ij}$ (in our setting, there is always at least one of the matrices that is symmetric).

The rest of this work is organized as follows: In Section \ref{sec:weakstrong}, which is the core of the paper, we prove Thereom \ref{thm:rentropy}. Section \ref{sec:logarithmicbounds} is dedicated to the proof of Theorem \ref{thm:gbounds}, which as stated above will consist in propagating (a stronger version of) the conditional logarithmic bounds required for Thereom \ref{thm:rentropy} to apply. It relies on classical parabolic arguments: a maximum principle and Schauder estimates, proofs of which are included in the appendix for the reader's convenience.

\subsection{Acknowledgments}
The author is thankful to Cyril Imbert and Clément Mouhot for many fruitful discussions on this work and related topics, and for proofreadings of this manuscript.

\section{The weak-strong relative entropy estimate}
\label{sec:weakstrong}
\subsection{Setup}

In this section we prove Theorem \ref{thm:rentropy}.
%We suppose $H(f_0\vert g_0)$ to be finite otherwise there is nothing to prove.
We divide the proof in three steps. The first one is Lemma \ref{lem:testlng}, which shows that we can use $\ln{g}$ as a test function in the weak formulation thanks to the bounds on $g$. We then compute the evolution of the relative entropy in Lemma \ref{lem:entropyevol}, isolating a negative good term and a bad error term. The good term is a relative version of the entropy production and can be used to control the worst part of the bad term. We conclude by a Gronwall argument.

We recall the propagation of moments for H-solutions (\cite{CarrapatosoDesvillettesHe2015}, Lemmas 7 and 8), which states that the moments of H-solutions of the Landau equation (and hence of classical solutions as well) grow at most linearly in time. We will only use that they are bounded on $[0,T]$, simplifying the result to, in our notations:
\begin{lemma}[Propagation of moments]
\label{lem:momentbounds}
    Consider the setting of Theorem \ref{thm:rentropy}. There exists two constants $M_{g}$ and $M_{f}$ both depending on $\gamma$, $T$ and $\rho$, and respectively on $M_{g_0}, H(g_0)$ and $M_{f_0},H({f_0})$, such that
    $$\sup_{t\in [0,T]} \int \jap{v}^{\rho-\gamma}g_t(v) dv \leq M_{g},$$
    $$\sup_{t\in [0,T]} \int \jap{v}^{\rho-\gamma}f_t(v) dv \leq M_{f}.$$
\end{lemma}
We will write our later results in terms of those constants rather than in terms of $M_{g_0}$ and $M_{f_0}$ to bound the moments of $f_t$ and $g_t$ we encounter.

To easily use the logarithmic bounds, we introduce non-negative functions $K_{gi}(t)$ for $i=1,2,3$, such that
\begin{align}
\label{eq:gradbound}
\vert \nabla \ln g_t\vert \leq K_{g1}(t) \jap{v}^{\kappa},\\
\label{eq:dtbound}
\vert \partial_t \ln g_t\vert \leq K_{g2}(t) \jap{v}^{\nu},\\
\label{eq:hessbound}
\Vert \nabla^2 \ln g_t\Vert \leq K_{g3}(t) \jap{v}^{\zeta},
\end{align}
and
$$\int_0^T (K_{g2}(t) + K_{g1}^2(t)+K_{g3}^2(t))dt <+\infty.$$
The finiteness of the above integral means that the $K_{gi}(t)$ are finite for almost every $t$. In practice, because of the regularizing effect of the Landau equation, the typical behavior of the $K_{gi}$ is to be decreasing from a possible blow up as $t\rightarrow0$.

One last point remains to be clarified before we begin the proof. In the following computations, we want to manipulate $\nabla f$, which is a priori a distribution, and $\nabla \ln f$ which is a priori not even defined. However, we know that H-solutions actually enjoy better regularity than what their definition seems to offer. We recall the entropy production estimate from \cite{Desvillettes2014a}:
\begin{lemma}[Entropy production estimate]
\label{lem:entropyprodestimate}
    For a probability density $f\in L^1$ with finite energy and entropy, there exists $c>0$ depending only on (upper bounds on) its energy and entropy, and on $\gamma$, such that
    $$\int_{\mathbb{R}^3} \vert\nabla \sqrt{f}\vert^2 \jap{v}^\gamma dv\leq c (D(f)+1).$$
\end{lemma}
This result was already mentioned in the introduction as the one used to show that H-solutions are weak solutions. Since $D(f_t)\in L^1([0,T])$ by hypothesis, it implies that $\nabla\sqrt{f}$ is a weighted $L^2([0,T],L^2)$ function, and not only a distribution. By a weighted Sobolev embedding, we also have that $\sqrt{f}$ lies in a weighted $L^2([0,T],L^6)$ space in velocity. This implies that the weak gradient $\nabla f = 2\nabla \sqrt{f} \sqrt{f}$ is a weighted $L^1([0,T],L^{3/2})$ function. In the H-formulation, we can hence simplify $(\nabla-\nabla')\sqrt{ff'}\sqrt{ff'}=\frac{1}{2}(\nabla-\nabla')ff'$ since it is now at least an $L^1_{loc}$ function.

However, $\nabla \ln f$ is still undefined and only $(\nabla \ln f) \sqrt{f}$ really makes sense. In the sequel, we will manipulate $\nabla \ln f$ quite freely because it makes the computations easier to follow, but one can check that in every integrand, any $\nabla \ln f$ can be factored with a $\sqrt{f}$ (and similarly for $\nabla \ln(ff')$ and $\sqrt{ff'}$) so that every term is well defined. The whole proof could in fact be rewritten without a single logarithm but would be much more cumbersome and cluttered.

\subsection{Using $\ln g$ as a test function}

In the spirit of many proofs of weak-strong uniqueness principles, we want to make use of the regularity of $g$ to use it as a test function in the weak formulation of the Landau equation. Since we work with the relative entropy, we are interested in the evolution of the quantity
$\int f \ln g $, so it is rather $\ln g$ that we need to be able to use as a test function. This is provided by the following lemma, which is essentially a technical application of the bounds on $\ln g$.

\begin{lemma}
\label{lem:testlng}
    In the setting of Theorem \ref{thm:rentropy}, $\ln g$ can be used as a test function in the weak formulation for $f$,
    \textit{i.e.} the following holds for all $t\in[0,T]$:
    \begin{align}
    \label{eq:testlng}
     \int_{\mathbb{R}^3} (\ln g_t)f_tdv &=\int_{\mathbb{R}^3} (\ln g_0)f_0dv + \int_0^t \int_{\mathbb{R}^3} (\partial_s\ln g_s)f_s dv\ ds \\
    & -\frac{1}{2}\int_0^t \iint_{\mathbb{R}^6} a(v-v') (\nabla\ln g_s - \nabla'\ln g'_s)\cdot\left((\nabla - \nabla')f_sf'_s\right) dvdv'ds.\nonumber
    \end{align}
    It is part of the result that every term is well-defined and finite.
\end{lemma}
\begin{proof}
\textit From the weak formulation, it is quite clear that it holds for test functions $\varphi$ that are only $C^1$ in time and $C^2$ in velocity (rather than $C^2$ in both). The smoothness and positivity of $g$ ensures that $\ln{g}$ satisfies this qualitative smoothness. We can thus build suitable test functions by truncating $\ln g$.

Consider a smooth bump function $\chi$, such that $0 \leq \chi \leq 1$, and which equals $1$ on $B(0,1)$ and vanishes outside $B(0,2)$. Fix a radius $R>0$, bound to go to infinity. The test function $$\varphi^R(t,v):= \chi\left(\frac{v}{R}\right)\ln g_t(v)$$ can be used in the weak formulation: it is smooth enough and compactly supported. Therefore \eqref{eq:testlng} holds with $\varphi^R$ instead of $\ln g$ and we now pass to the limit in each term.
\\

We do a little trick to show convergence and finiteness for \textit{every} $t$ of $\int \ln g_t f_t$ even though we have only $L^1$ in time control. Pick a $t_0$ such that $K_{g1}(t_0)$ is finite. Integrating the gradient bound \eqref{eq:gradbound} on $\ln g$ from $0$ to $v$ yields
\begin{equation}
\label{eq:boundongfromgrad}
\vert \ln g_{t_0}(v)\vert \leq \vert \ln g_{t_0}(0)\vert + K_{g1}(t_0) \jap{v}^{\kappa}\vert v \vert \leq C_g \jap{v}^{\kappa+1}.
\end{equation}
Then using the bound \eqref{eq:dtbound} on the time derivative,
$$\vert \ln g_{t}(v)\vert \leq \vert \ln g_{t_0}(v)\vert+\vert t-t_0 \vert \int_{t_0}^{t} K_{g2}(s)ds \jap{v}^\nu \leq C_g \jap{v}^{\max(\kappa+1,\nu)}.$$
The pointwise convergence $\varphi^R \rightarrow \ln g$ as $R \rightarrow \infty$ and the bound $\vert\varphi^R\vert \leq C_{g}\jap{v}^{\max(\kappa+1,\nu)}$ together with the control of the $\max(\kappa+1,\nu)$-th moment of $f$ ensures that, by dominated convergence, for all $t\in [0,T]$,
$$\int_{\mathbb{R}^3} \varphi^R_t f_t dv \xrightarrow[R\rightarrow \infty]{} \int_{\mathbb{R}^3} \ln(g_t) f_t dv,$$%\in [-\infty,\ln \Vert g\Vert_\infty]$$
with a finite right hand side.
%For $t=0$, the finiteness of $H(f_0)$ and $H(f_0\vert g_0)$ implies that $\int \ln g_0 f_0$ is finite.
\\

%For the other terms we use dominated convergence.
Similarly, the logarithmic time derivative bound $\vert \partial_s\varphi^R - \partial_s\ln g \vert \leq K_{g2}(s) \jap{v}^\nu$ (true thanks to \eqref{eq:dtbound}) and the uniform-in-time control of the $\nu$-th moment of $f$ ensures that
$$\int_0^t\int_{\mathbb{R}^3} (\partial_s\varphi^R_s) f_s dvds \xrightarrow[R\rightarrow \infty]{} \int_0^t\int_{\mathbb{R}^3} (\partial_s\ln g_s) f_s dvds$$
because the domination $K_{g2}(s)\jap{v}^\nu f_s$ is integrable on $[0,t]\times \mathbb{R}^3$.
\\

It remains to pass to the limit in the most complicated term, that is
$$\int_0^t \iint_{\mathbb{R}^6} a(v-v') (\nabla\varphi^R_s - \nabla'(\varphi^R_s)')\cdot(\nabla - \nabla')f_sf'_s\ dvdv'ds.$$
We rely on the bound in Remark \ref{rem:weaksense}, which makes it enough to show that
$$\int_0^t\iint_{\mathbb{R}^6}  f_sf'_s \vert\nabla\psi^R_s - \nabla'(\psi^R_s)'\vert^2\vert v-v'\vert^{\gamma +2}dvdv'ds \xrightarrow[R\rightarrow \infty]{} 0,$$
where $\psi^R = \ln g -\varphi^R$. We first treat the easier case $\gamma+2 \geq 0$. In this case, there is no singularity. We have 
$$\nabla \psi^R =\left(1-\chi\left(\frac{\cdot}{R}\right)\right)\nabla \ln g  - \frac{1}{R}\nabla \chi\left(\frac{\cdot}{R}\right) \ln g.$$
Using the bound \eqref{eq:gradbound} for $\nabla \ln g$ and \eqref{eq:boundongfromgrad} for $\ln g$, and the support of $\chi$ we get
$$\vert \nabla \psi^R_s \vert \leq C K_{g1}(s) \left( \jap{v}^\kappa  + \frac{\jap{v}^{\kappa+1}}{R}\mathbf{1}_{\vert v \vert \leq 2R}\right)\leq C K_{g1}(s) \jap{v}^\kappa.$$
Hence
\begin{align*}
    \vert\nabla\psi_s^R - \nabla'(\psi_s^R)'\vert^2\vert v-v'\vert^{\gamma +2} &\leq C K_{g1}^2(s) \left( \jap{v}^{2\kappa+\gamma+2}  + \jap{v'}^{2\kappa+\gamma+2}\right)\\
    &\leq C K_{g1}^2(s) \left( \jap{v}^{\rho}  + \jap{v'}^{\rho}\right),
\end{align*}
since $2\kappa+\gamma+2 \leq 2\kappa+2\gamma+4\leq \rho.$
As before, because of the control of the $\rho$-th moment of $f$ and by hypothesis on $K_{g1}(s)$, this last term is integrable on $[0,t]\times \mathbb{R}^3$, so the pointwise convergence $\nabla \psi^R\rightarrow 0$ concludes.

We now treat the case $\gamma+2 <0$. Both the tail and the singularity $\{v=v'\}$ must be dealt with.
It is easy to see that the non-singular part goes to zero, \textit{i.e.}
$$\int_0^t\iint_{\vert v-v' \vert \geq 1}  f_sf'_s \vert\nabla\psi^R_s - \nabla'(\psi^R_s)'\vert^2\vert v-v'\vert^{\gamma +2}dvdv'ds \xrightarrow[R\rightarrow \infty]{} 0,$$
because of the bounds
$$ \vert\nabla\psi^R_s - \nabla'(\psi^R_s)'\vert^2 \leq C K_{g1}^2(s)\left( \jap{v}^{2\kappa}  + \jap{v'}^{2\kappa}\right) $$
and $\vert v-v'\vert^{\gamma +2}\leq 1$, and the control of the moments of $f$.

For the singular part $\{\vert v - v'\vert <1 \}$, we compute the Hessian of $\psi^R$:
\begin{align*}
    \nabla^2 \psi^R_s = \left(1-\chi\left(\frac{\cdot}{R}\right)\right)\nabla^2 \ln g  - \frac{2}{R}\nabla \chi\left(\frac{\cdot}{R}\right) \otimes \nabla\ln g -\frac{1}{R^2}\nabla^2 \chi\left(\frac{\cdot}{R}\right) \ln g.
\end{align*}
We use the growth bounds on $g$ and its derivatives (namely, \eqref{eq:gradbound}, \eqref{eq:hessbound}  and \eqref{eq:boundongfromgrad}), to get
\begin{align*}
\Vert \nabla^2 \psi^R_s \Vert &\leq C(K_{g1}(s)+K_{g3}(s)) \left(\jap{v}^\zeta  + \frac{\jap{v}^\kappa}{R}\mathbf{1}_{\vert v \vert \leq 2R}+\frac{\jap{v}^{\kappa+1}}{R^2}\mathbf{1}_{\vert v \vert \leq 2R} \right) \\
&\leq C(K_{g1}(s)+K_{g3}(s))\jap{v}^\frac{\rho}{2}.
\end{align*}
From this we obtain:
$$\vert\nabla\psi^R_s - \nabla'(\psi^R_s)'\vert \leq C(K_{g1}(s)+K_{g3}(s)) \vert v-v'\vert \left( \jap{v}^\frac{\rho}{2} + \jap{v'}^\frac{\rho}{2}\right), $$
a bound we can use to cancel the singularity: if $\vert v - v'\vert < 1$, then
\begin{align*}
    \vert\nabla\psi^R - \nabla'(\psi^R)'\vert^2\vert v-v'\vert^{\gamma +2} &\leq C(K_{g1}^2(s)+K_{g3}^2(s)) \vert v-v'\vert^{\gamma +4}  \left( \jap{v}^{\rho}  + \jap{v'}^{\rho}\right)\\ &\leq C(K_{g1}^2(s)+K_{g3}^2(s))\left( \jap{v}^{\rho}  + \jap{v'}^{\rho}\right),
\end{align*} 
because $\gamma+4>0$. We conclude as above thanks to the control of the $\rho$-th moment of $f$. Hence the weak formulation holds with test function $\ln g$.% The finiteness of every other term implies that $\int \ln(g_t) f_t$ is finite.
\end{proof}
% \begin{remark}
%     The fact that $\int \ln(g_t) f_t$ is finite can be more directly seen as a consequence of the estimate on $\vert \nabla \ln g_t\vert$, which in turn implies a polynomial bound on $\vert \ln g_t \vert$ that makes it integrable against $f_t$.
% \end{remark}

\subsection{Evolution of the relative entropy}
We now proceed to the main computation. We write the evolution of entropy as the sum of a bad error term minus a positive good term.
\begin{lemma}
    \label{lem:entropyevol}
Consider the setting of Theorem \ref{thm:rentropy}. The evolution of the relative entropy satisfies:
\begin{equation}
    \label{eq:entropyevol}
    H(f_t \vert g_t) \leq H(f_0 \vert g_0) + \int_0^t (-\mathcal{G}_s+\mathcal{B}_s) ds
\end{equation}
with the good term being
\begin{equation*}
    \mathcal{G}_s:=\frac{1}{2}\iint a(v-v'):\left[(\nabla-\nabla')\ln\left(\frac{f_sf'_s}{g_sg'_s}\right)\right]^{\otimes2}f_sf'_s dvdv'\geq 0
\end{equation*}
and the bad term
\begin{equation*}
    \mathcal{B}_s:=\frac{1}{2}\iint \frac{1}{g_sg_s'}(f_s-g_s)(f_s'-g_s')(\nabla-\nabla') \cdot \left[ a (v - v') (\nabla - \nabla')(g_sg_s')\right]dvdv'.
\end{equation*}
\end{lemma}
\begin{remark}
     We recall that $\nabla \ln f$ has a sense only when paired with $\sqrt{f}$, and is then essentially a shorthand notation for $2\nabla \sqrt{f}$. For instance, the above term $\mathcal{G}_s$ rigorously means
 \begin{align*}
    \mathcal{G}_s&=2\iint a(v-v'):\left[(\nabla-\nabla')\sqrt{f_sf_s'}-\sqrt{f_sf_s'}(\nabla-\nabla')\ln(g_sg'_s) \right]^{\otimes2} dvdv'.
\end{align*}    
\end{remark}
\begin{proof}
Thanks to Lemma \ref{lem:testlng}, we can write
$$H(f_t \vert g_t) = H(f_t) - \int_{\mathbb{R}^3} f_t \ln g_tdv$$
and use the weak formulation \eqref{eq:testlng} for the second term. Also plugging in the entropy production estimate for the first term,
we have:
\begin{align*}
    H(f_t \vert g_t) &\leq  H(f_0)-\frac{1}{2}\int_0^t\iint_{\mathbb{R}^6}a (v - v'): \left[(\nabla - \nabla')\ln(f_sf'_s) \right]^{\otimes 2}f_sf'_sdvdv'ds \\
    &-\int_{\mathbb{R}^3} f_0 \ln g_0 - \int_0^t \int_{\mathbb{R}^3} (\partial_s\ln g_s)f_s dv\ ds \\
    & +\frac{1}{2}\int_0^t \iint_{\mathbb{R}^6} a(v-v') (\nabla\ln g_s - \nabla'\ln g'_s)\cdot(\nabla - \nabla')f_sf'_s\ dvdv'ds \\
    &\leq H(f_0 \vert g_0) + \int_0^t \mathcal{F}_s ds
\end{align*}
where
\begin{align*}
    \mathcal{F}_s:=&-\frac{1}{2}\iint_{\mathbb{R}^6}a (v - v') :\left[(\nabla - \nabla')\ln(f_sf'_s) \right]^{\otimes 2}f_sf'_sdvdv' - \int_{\mathbb{R}^3} (\partial_s\ln g_s)f_s dv \\
    & +\frac{1}{2} \iint_{\mathbb{R}^6} a(v-v') (\nabla\ln g_s - \nabla'\ln g'_s)\cdot (\nabla - \nabla')f_sf'_s\ dvdv'.
\end{align*}
We drop the time index $s$ from now on and want to show that $\mathcal{F}=-\mathcal{G}+\mathcal{B}$. Using the identities
$$(\nabla - \nabla')ff' = \left((\nabla - \nabla')\ln(ff')\right) ff'$$
and
$$\nabla\ln g - \nabla'\ln g' = (\nabla - \nabla')\ln(gg'),$$
we can combine terms appropriately in order to write:
\begin{align*}
    \mathcal{F}=&-\frac{1}{2}\iint_{\mathbb{R}^6}\left((\nabla - \nabla')\ln\left(\frac{ff'}{gg'}\right)\right)\cdot a (v - v') \left((\nabla - \nabla')\ln(ff')\right) ff'dvdv' \\
    &- \int_{\mathbb{R}^3} (\partial_s\ln g)\ f dv.
\end{align*}
Let us introduce a notation for the semidefinite quantity
$$\Phi(f,g) :=\left((\nabla-\nabla')\ln(ff') \right) \cdot a (v - v') (\nabla - \nabla')\ln(gg'),$$
so that we can recast $\mathcal{F}$ in the more compact form:
\begin{align}
\label{eq:F}
    \mathcal{F}=&-\frac{1}{2}\iint_{\mathbb{R}^6}\Phi\left(\frac{f}{g}, f\right) ff'dvdv'- \int_{\mathbb{R}^3} (\partial_s\ln g) f dv.
\end{align}
The first term we keep as is and we now work on the second term: since $\int \partial_s g=0$ because it is a divergence, we can write
\begin{align*}
    \int_{\mathbb{R}^3} (\partial_s\ln g) f dv &=  \int_{\mathbb{R}^3} \frac{\partial_s g}{g}\ f dv\\
    &= \int_{\mathbb{R}^3} \frac{\partial_s g}{g}\ (f-g) dv.
\end{align*}
Using that $g$ is a classical solution of the Landau equation,
\begin{align*}
    \int_{\mathbb{R}^3} \frac{\partial_s g}{g}\ (f-g) dv&=
    \iint \frac{1}{g}(f-g)(\nabla-\nabla') \cdot \left[ a (v - v') (\nabla - \nabla')(gg')\right]dvdv',
\end{align*}
where we have added an artificial $\nabla'\cdot$ which cancels out in the $dv'$ integration. We wish to symmetrize the integrand. Consider the following identity:
\begin{align*}
    \frac{1}{g}(f-g)
    &=\frac{1}{2gg'}(2fg'-gg'-ff') + \frac{1}{2gg'}(ff'-gg').
\end{align*}
Only the first term is not symmetric in primed and non-primed quantities. This means that the symmetrized version is
\begin{align}
\label{eq:symmetrized}
    \frac{1}{2} \left[\frac{1}{g}(f-g)  + \frac{1}{g'}(f'-g')\right]&=\frac{1}{2gg'}(fg'+f'g-gg'-ff') + \frac{1}{2gg'}(ff'-gg')\nonumber\\
    &=-\frac{1}{2gg'}(f-g)(f'-g') + \frac{1}{2gg'}(ff'-gg').
\end{align}
Plugging this above, we get that
\begin{align*}
    \int_{\mathbb{R}^3} \partial_s\ln g\ f dv =&-
    \frac{1}{2}\iint \frac{1}{gg'}(f-g)(f'-g')(\nabla-\nabla') \cdot \left[ a (v - v') (\nabla - \nabla')(gg')\right]dvdv'\\
    &+\frac{1}{2}\iint \frac{1}{gg'}(ff'-gg')(\nabla-\nabla') \cdot \left[ a (v - v') (\nabla - \nabla')(gg')\right]dvdv'\\
    =&-\mathcal{B}+\frac{1}{2}\iint \frac{1}{gg'}(ff'-gg')(\nabla-\nabla') \cdot \left[ a (v - v') (\nabla - \nabla')(gg')\right]dvdv'.
\end{align*}
The second term of this expression can be integrated by parts:
\begin{align*}
    \frac{1}{2}\iint \frac{1}{gg'}(ff'-gg')&(\nabla-\nabla') \cdot \left[ a (v - v') (\nabla - \nabla')(gg')\right]dvdv'\\
    &=-\frac{1}{2}\iint (\nabla-\nabla')\left[\frac{1}{gg'}(ff'-gg')\right]  \cdot a (v - v') (\nabla - \nabla')(gg')dvdv'.
\end{align*}
Using the identities
$$(\nabla-\nabla')\left[\frac{1}{gg'}(ff'-gg')\right] =\left[(\nabla-\nabla') \ln\left(\frac{ff'}{gg'}\right)\right]\frac{ff'}{gg'} $$
and
$$(\nabla - \nabla')(gg')=\left((\nabla - \nabla')\ln(gg')\right) gg',$$
we can express this integral using $\Phi$:
\begin{align*}
   -\frac{1}{2}\iint (\nabla-\nabla')&\left[\frac{1}{gg'}(ff'-gg')\right]\cdot  a (v - v') (\nabla - \nabla')(gg')dvdv'\\
   &=-\frac{1}{2}\iint\left( (\nabla-\nabla') \ln\left(\frac{ff'}{gg'}\right)\right)\cdot a (v - v') \left((\nabla - \nabla')\ln(gg')\right)ff'dvdv'\\
   &= -\frac{1}{2}\iint \Phi\left(\frac{f}{g},g\right) ff'dvdv'.
\end{align*}
To summarize, we have reached the expression
\begin{align*}
    \int_{\mathbb{R}^3} (\partial_s\ln g)\ f dv =-\mathcal{B}-\frac{1}{2}\iint \Phi\left(\frac{f}{g},g\right) ff'dvdv'.
\end{align*}
Going back to the equation \eqref{eq:F} for $\mathcal{F}$,
\begin{align*}
    \mathcal{F}&=-\frac{1}{2}\iint_{\mathbb{R}^6}\Phi\left(\frac{f}{g}, f\right) ff'dvdv'+\frac{1}{2}\iint \Phi\left(\frac{f}{g},g\right) ff'dvdv'+\mathcal{B}\\
    &=-\frac{1}{2}\iint_{\mathbb{R}^6}\Phi\left(\frac{f}{g}, \frac{f}{g}\right) ff'dvdv'+\mathcal{B}\\
    &= -\mathcal{G} + \mathcal{B}.
\end{align*}
This is the desired expression.
\end{proof}
\begin{remark}
    The good term $\mathcal{G}$ was expected to appear: it is a relative version of the entropy production. It was formed by combining the other terms with the nicely tensorized $(ff'-gg')$ part of the right hand side of \eqref{eq:symmetrized}. The bad term $\mathcal{B}$ comes from the remaining term $(f-g)(f'-g')$ in \eqref{eq:symmetrized}, so it can really be interpreted as a "tensorization defect".
\end{remark}
\subsection{Control of the bad term by the good one}

We first state what we will refer to as the \textit{refined Pinsker's inequality}:
\begin{lemma}
\label{lem:pinsker}
    For all probabilities $f,g$ with densities,
    $$\int \left\vert\sqrt{f}-\sqrt{g}\right\vert^2 \leq H(f\vert g).$$
\end{lemma}
\begin{proof}
Making use of $-\ln x\geq 1-x$,
\begin{align*}
    H(f\vert g) &= -2\int f\ln\left(\sqrt{\frac{g}{f}}\right)\\
    &\geq 2\int f\left(1-\sqrt{\frac{g}{f}}\right)\\
    &=2 -\int 2\sqrt{fg}\\
    &=\int f +g-2\sqrt{fg}\\
    &=\int \left\vert\sqrt{f}-\sqrt{g}\right\vert^2,
\end{align*}
using in lines 3 and 4 that $f$ and $g$ have unit mass.
\end{proof}
The good term $-\mathcal{G}$ is non-positive. However we can do better and show that, up to an error, it controls a weighted relative Fisher information. This is due to the following key classical fact: the convolution of the non-negative matrix $a$ with $f$ yields a \textit{positive} matrix, because of the finite entropy of $f$. This holds thanks to the behavior of $a$ at infinity and not because of its singularity, so it remains true for a regularized cut-off version of $a$, as stated in the following lemma:
\begin{lemma}
\label{lem:ellipticitycutoff}
    Consider the cut-off kernel $$\tilde{a}(z):=\jap{z}^\gamma \vert z \vert^2 \Pi(z).$$
    It is smooth on $\mathbb{R}^3$, satisfies $a(z)\geq \tilde{a}(z)$ and the following ellipticity estimate: for any $s\in[0,T]$
    $$\tilde{a}\star f_s(v):=\int\tilde{a}(v-v')f_s(v')dv' \geq c_0 \jap{v}^\gamma \Id$$
    where $\Id$ is the $3\times 3$ identity matrix, and $c_0$ depends only on (upper bounds on) the energy and entropy of $f_0$.
\end{lemma}
The proof of this ellipticity estimate and the fact that it holds even after cut-off (with the same proof) is classical (see for instance \cite{GolseImbertJi2024}, Lemma 2.1).

Using this cut off kernel, we are able to exploit the coercivity of $\mathcal{G}$ with a well-behaved error term. The next lemma can be thought of as a relative version of Lemma \ref{lem:entropyprodestimate}.
\begin{lemma}
    \label{lem:good}
We keep the setting of Theorem \ref{thm:rentropy}, and recall that, omitting time indices,
$$\mathcal{G}=\frac{1}{2}\iint a(v-v'):\left[(\nabla-\nabla')\ln\left(\frac{ff'}{gg'}\right)\right]^{\otimes2}ff' dvdv'.$$
It satisfies
\begin{align}
    \label{eq:good}
    \mathcal{G} \geq c_0\int \left\vert\nabla\ln\left(\frac{f}{g} \right) \right\vert^2 \langle v \rangle^{\gamma}f dv\ - C(1+K_{g2}^2)(M_{f}+M_{g})H(f\vert g)
\end{align}
where $c_0$ depends only on the energy and entropy of $f_0$, and $C$ depends only on $\gamma$.
\end{lemma}
\begin{proof}
  We use the cut-off kernel $\tilde{a}(v)$ from Lemma \ref{lem:ellipticitycutoff} which is smooth (even at the origin), and behaves like $a(v)$ at infinity. Using that $a(v-v')\geq \tilde{a}(v-v')$, we bound from below the good term. Following a classical technique of the non relative setting (see for instance \cite[Section 2]{Ji2024b} or \cite[Proposition 2.2]{GolseImbertJi2024}) we cut $\mathcal{G}$ in a coercive part and an error:
\begin{align*}
    \mathcal{G}
    &\geq\frac{1}{2}\iint \tilde{a}(v-v'):\left[(\nabla-\nabla')\ln\left(\frac{ff'}{gg'}\right)\right]^{\otimes2}ff' dvdv'\\
    &=\frac{1}{2}\iint \tilde{a}(v-v'):\left[\nabla\ln\left(\frac{f}{g} \right)-\nabla'\ln\left(\frac{f'}{g'} \right)\right]^{\otimes2}ff'dvdv'\\
    &=\iint \tilde{a}(v-v'): \left[\nabla\ln\left(\frac{f}{g} \right)\right]^{\otimes2}ff'dvdv' - \iint \tilde{a}(v-v') \nabla\ln\left(\frac{f}{g} \right)\cdot \nabla'\ln\left(\frac{f'}{g'} \right)ff'dvdv'
\end{align*}
by symmetry.
The integral in $dv'$ in the first term can be rewritten in terms of the convolution $\tilde{a}\star f$. We integrate by parts the second term (in both $v$ and $v'$).
\begin{align*}
   \mathcal{G}&\geq\int \tilde{a}\star f : \left[\nabla\ln\left(\frac{f}{g} \right)\right]^{\otimes2}fdv - \iint \tilde{a}(v-v') \nabla\left(\frac{f}{g} \right)\cdot \nabla'\left(\frac{f'}{g'} \right)gg'dvdv'\\
    &\geq c_0\int \left\vert\nabla\ln\left(\frac{f}{g} \right) \right\vert^2 \langle v \rangle^{\gamma}f dv\ -\iint \nabla' \cdot \nabla \cdot (\tilde{a}(v-v')gg') \left(\frac{f}{g}-1 \right)\left(\frac{f'}{g'} -1\right)dvdv'
\end{align*}
The first term was obtained by applying Lemma \ref{lem:ellipticitycutoff} and is the weighted relative Fisher information we seek. The second one is the error term and we wish to control it by the relative entropy.
We compute
\begin{align*}
\nabla' \cdot \nabla \cdot (\tilde{a}(v-v')g(v)g(v')) &=-(\nabla\cdot \nabla \cdot \tilde{a})(v-v')g(v)g(v')+(\nabla \cdot\tilde{a})(v-v')g(v)\nabla g(v')\\
&\ \ \ \ \ - (\nabla \cdot\tilde{a})(v-v')\nabla g(v) g(v') + \tilde{a}(v-v')\nabla g(v)\nabla g(v').
\end{align*}
A direct computation on the smooth $\tilde{a}$ yields the bounds:
$$ \Vert \tilde{a} \Vert \leq C\jap{v}^{\gamma+2},\; \; \;  \; \; \;
\vert \nabla\cdot \tilde{a} \vert \leq C\jap{v}^{\gamma+1},\; \; \; 
\; \; \;  \vert \nabla\cdot\nabla\cdot\tilde{a} \vert \leq C\jap{v}^{\gamma}.$$
Note that $C$ depends only on $\gamma$. Thanks to these bounds, we get the estimate
\begin{align*}
    \vert \nabla' \cdot \nabla \cdot (\bar{a}gg') \vert \leq &C\bigg(\jap{v-v'}^{\gamma}+\vert \nabla \ln(g)\vert\jap{v-v'}^{\gamma+1}\\
    &+\vert \nabla' \ln(g')\vert\jap{v-v'}^{\gamma+1}+\vert \nabla \ln(g)\vert\vert \nabla '\ln(g')\vert\jap{v-v'}^{\gamma+2}\bigg) gg'.
\end{align*}
Denoting $\alpha_+=\max(\alpha,0)$ for a real number $\alpha$, it always hold that
$$\jap{v-v'}^\alpha \leq \jap{v}^{\alpha_+}+\jap{v'}^{\alpha_+}.$$
Plugging this in as well as the gradient growth bound $\vert \nabla \ln(g)\vert \leq K_{g1}\jap{v}^\kappa$, we end up with
\begin{align*}
    \vert \nabla' \cdot \nabla \cdot (\bar{a}gg') \vert \leq &C\left(1+K_{g1}^2\right)\bigg(\jap{v}^{\gamma_+}+\jap{v'}^{\gamma_+}+\left(\jap{v}^\kappa+\jap{v'}^\kappa\right)\left(\jap{v}^{(\gamma+1)_+}+\jap{v'}^{(\gamma+1)_+}\right)\\
    &+\left(\jap{v}^\kappa\jap{v'}^\kappa\right)\left(\jap{v}^{(\gamma+2)_+}+\jap{v'}^{(\gamma+2)_+}\right)\bigg)gg'\\
    \leq & C\left(1+K_{g1}^2\right)\ \jap{v}^{\kappa+(\gamma+2)_+}\ \jap{v'}^{\kappa+(\gamma+2)_+} gg'\\
    \leq &C\left(1+K_{g1}^2\right)\jap{v}^{\frac{\rho}{2}}\jap{v'}^{\frac{\rho}{2}} gg'.
\end{align*}
Plugging this estimate in the error term, the resulting integrand factors out in primed and unprimed quantities:
\begin{align*}
    \iint \vert \nabla' \cdot \nabla \cdot (\bar{a}gg') \vert \left\vert\frac{f}{g}-1 \right\vert&\left\vert\frac{f'}{g'} -1\right\vert dv dv'\\
    &\leq C\left(1+K_{g1}^2\right)\iint \jap{v}^{\frac{\rho}{2}}\jap{v'}^{\frac{\rho}{2}}  \vert f-g \vert \vert f'-g'\vert dvdv' \\
    &\leq C\left(1+K_{g1}^2\right) \left(\int \jap{v}^{\frac{\rho}{2}} \vert f-g \vert dv\right)^2.
\end{align*}
Writing $f-g=(\sqrt{f}-\sqrt{g})(\sqrt{f}+\sqrt{g})$ and using the Cauchy-Scwharz inequality,
\begin{align*}
    \left(\int \jap{v}^{\frac{\rho}{2}} \vert f-g \vert dv\right)^2 &\leq \left(\int \jap{v}^{\rho}\vert \sqrt{f}+\sqrt{g}\vert ^2 dv\right) \left( \int \vert \sqrt{f}-\sqrt{g}\vert ^2dv \right)\\
    &\leq 2 \left(\int \jap{v}^{\rho}(f+g)dv\right) H(f\vert g)\\
    &\leq 2 (M_f +M_g) H(f\vert g),
\end{align*}
using the refined Pinsker's inequality (Lemma \ref{lem:pinsker}) and the control of the $\rho$-th moment of $f$ and of $g$.

Gathering everything, the good term satisfies:
\begin{align*}
    \label{eq:good}
    \mathcal{G} \geq c_0\int \left\vert\nabla\ln\left(\frac{f}{g} \right) \right\vert^2 \langle v \rangle^{\gamma}f dv\ - C(1+K_{g1}^2)(M_{f}+M_{g})H(f\vert g),
\end{align*}
which is the desired estimate.
\end{proof}

\begin{remark}
    The cut-off $\tilde{a}$ was introduced so that the error term is reasonable: without it, the double divergence $\nabla\cdot\nabla\cdot a$ would appear, which is a multiple of a Dirac mass when $\gamma=-3$. The error term would then be an $L^2$ norm of $(f-g)$ which we cannot control by the relative entropy.
\end{remark}

We now analyse the bad term, making the amount of relative Fisher information that we control by $\mathcal{G}$ appear:
\begin{lemma}
    \label{lem:bad}
We keep the setting of Theorem \ref{thm:rentropy}, and recall that, omitting time indices,
$$\mathcal{B}=\frac{1}{2} \iint \frac{1}{gg'}( f-g ) ( f'-g') (\nabla-\nabla')\cdot \left[ a (v - v') (\nabla - \nabla')(gg')\right]dvdv'.$$
It satisfies
\begin{align}
    \label{eq:bad}
    \vert\mathcal{B}\vert \leq c_0\int \left\vert\nabla \ln \left(\frac{f}{g}\right) \right\vert^2 \langle v \rangle^{\gamma}fdv + C\frac{(K_{g1}^2+K_{g3}^2) (M^2_f+M^2_{g})}{c_0}H(f\vert g),
\end{align}
with the same $c_0$ as in Lemma \ref{lem:good}, and $C$ depending only on $\gamma$.
\end{lemma}
\begin{proof}
To make the relative Fisher information appear, we un-symmetrize the gradient, then integrate by parts in $v$, and then write everything with logarithmic derivatives:
\begin{align*}
    \mathcal{B} &= \iint \left( \frac{f}{g}-1 \right) \left( \frac{f'}{g'}-1 \right) \nabla \cdot \left[ a (v - v') (\nabla - \nabla')(gg')\right]dvdv'\\
    &= -\iint \left( \frac{f'}{g'} - 1\right)\nabla \left( \frac{f}{g} - 1\right)\cdot a (v - v') (\nabla - \nabla')(gg')dvdv'\\
    &=-\iint f\left( f'-g'\right)\nabla \ln \left(\frac{f}{g}\right) \cdot a (v - v') (\nabla - \nabla')\ln(gg')dvdv'.\\
\end{align*}
Using $f'-g'=(\sqrt{f'}-\sqrt{g'})(\sqrt{f'}+\sqrt{g'})$, by Cauchy-Schwarz and introducing the right weights, we have:
\begin{align}
\label{eq:expressionb}
    \vert \mathcal{B} \vert &\leq \iint \left\vert\nabla \ln \left(\frac{f}{g}\right) \right\vert f\left\vert \sqrt{f'}-\sqrt{g'}\right\vert\left\vert \sqrt{f'}+\sqrt{g'}\right\vert \left\vert a (v - v') (\nabla - \nabla')\ln(gg') \right\vert dvdv'\nonumber\\
    &\leq \left(\int \left\vert\nabla \ln \left(\frac{f}{g}\right) \right\vert^2 \langle v \rangle^{\gamma}fdv\right)^\frac{1}{2}\times \left(\int \left\vert \sqrt{f'}-\sqrt{g'}\right\vert^2dv'\right)^\frac{1}{2}\nonumber\\
    &\;\;\;\;\; \times \left(\iint f\langle v \rangle^{-\gamma}\left\vert\sqrt{f'}+\sqrt{g'}\right\vert^2\left\vert a (v - v') (\nabla - \nabla')\ln(gg') \right\vert^2dvdv'\right)^\frac{1}{2}.
\end{align}
Notice that we were again able to split one of the double integrals since its integrand was a product of primed and non-primed quantities.
The first integral is the relative Fisher information we sought, the second integral will be controlled thanks to the refined Pinsker's inequality. To deal with the last integral, we proceed in a similar fashion as back in Lemma \ref{lem:testlng}: we want to use the $C^1$ regularity of $\nabla \ln g$ to cancel out the potential singularity of $a(v-v')=\vert v-v'\vert^{\gamma+2}\Pi(v-v')$.

We first suppose that $\gamma+2<0$. We make use of the Hessian estimate \eqref{eq:hessbound} to write
$$\vert\nabla\ln g - \nabla'\ln g'\vert \leq C K_{g3} \vert v-v'\vert \left( \jap{v}^\zeta + \jap{v'}^\zeta\right). $$
This yields
$$\vert a (v - v') (\nabla - \nabla')\ln(gg') \vert \leq C K_{g3} \vert v-v'\vert^{\gamma+3} \left( \jap{v}^\zeta + \jap{v'}^\zeta\right),$$
which we use for the singular part $\vert v - v' \vert<1$ (over which $\vert v-v'\vert^{\gamma+3}\leq 1$). 
%We note that this is the only point in the whole proof where we use that $\gamma\geq -3$, everything else works with $\gamma>-4$.
If $\vert v - v' \vert \geq 1$, we simply use the gradient estimate $\vert \nabla \ln g \vert\leq K_{g1} \jap{v}^{\kappa}$ and we can drop the $\vert v-v'\vert^{\gamma+2}$ because it is less than $1$, yielding
$$\vert a (v - v') (\nabla - \nabla')\ln(gg') \vert \leq C K_{g1} (\jap{v}^{\kappa}+\jap{v'}^{\kappa})$$
for $\vert v - v' \vert \geq 1$.

If $\gamma+2\geq 0$, there is no singularity and the estimate
$$\vert a (v - v') (\nabla - \nabla')\ln(gg') \vert \leq C K_{g1} (\jap{v}^{\kappa+\gamma+2}+\jap{v'}^{\kappa+\gamma+2})$$
holds. 

In both cases, it always holds for all $v,v'$ that
$$\vert a (v - v') (\nabla - \nabla')\ln(gg')\vert \leq C (K_{g1}+K_{g3}) \left( \jap{v}^\frac{\rho}{2} + \jap{v'}^\frac{\rho}{2}\right). $$
Hence, the last integral in the previous expression \eqref{eq:expressionb} for $\mathcal{B}$ can be bounded using moments of $f$ and $g$:
\begin{align*}
 \iint f\langle v \rangle^{-\gamma}\left\vert\sqrt{f'}+\sqrt{g'}\right\vert^2&\left\vert a (v - v') (\nabla - \nabla')\ln(gg') \right\vert^2dvdv'\\
 &\leq  C (K_{g1}+K_{g3})^2 \iint f(f'+g') \jap{v}^{-\gamma}\left( \jap{v}^{\rho} + \jap{v'}^{\rho}\right)dvdv'\\
 &\leq C (K_{g1}+K_{g3})^2 (M_f+M_{g})^2.
\end{align*}
After plugging this in \eqref{eq:expressionb}, we use Young's inequality on the two remaining integrals: $2a^\frac{1}{2}b^\frac{1}{2}\leq \epsilon a +\frac{1}{\epsilon}b$ for any $\epsilon>0$. It allows us to write:
\begin{align*}
    \vert \mathcal{B} \vert 
    &\leq C(K_{g1}+K_{g3}) (M_f+M_{g})\left(\int \left\vert\nabla \ln \left(\frac{f}{g}\right) \right\vert^2 \langle v \rangle^{\gamma}fdv\right)^\frac{1}{2}\left(\int \left\vert \sqrt{f'}-\sqrt{g'}\right\vert^2dv'\right)^\frac{1}{2}\\
    &\leq c_0\int \left\vert\nabla \ln \left(\frac{f}{g}\right) \right\vert^2 \langle v \rangle^{\gamma}fdv + C\frac{(K_{g1}^2+K_{g3}^2) (M^2_f+M^2_{g})}{c_0}\int \left\vert \sqrt{f'}-\sqrt{g'}\right\vert^2dv'\\
    &\leq c_0\int \left\vert\nabla \ln \left(\frac{f}{g}\right) \right\vert^2 \langle v \rangle^{\gamma}fdv + C\frac{(K_{g1}^2+K_{g3}^2) (M^2_f+M^2_{g})}{c_0}H(f\vert g),
\end{align*}
where in the last line we used the refined Pinsker's inequality (Lemma \ref{lem:pinsker}). This is the bound we wanted.
\end{proof}

We are now all set to state and prove a precise form of the estimate in Theorem \ref{thm:rentropy}.
\begin{prop}
\label{prop:gronwall}
    Consider the setting of Theorem \ref{thm:rentropy}.
    The relative entropy estimate \eqref{eq:thm_rentropy} holds. More precisely, for any $t\in[0,T]$,
    $$H(f_t \vert g_t) \leq H(f_0 \vert g_0) \exp\left(\mathcal{C}\int_0^t (K_{g1}^2(s)+K_{g3}^2(s))ds\right)$$
    with
    $$\mathcal{C}=C_\gamma(1+M_{f}^2+M_{g}^2)(1+c_0^{-1})$$
    where $c_0^{-1}$ depends on the energy and entropy of $f_0$.
\end{prop}

\begin{proof} 

Recall Lemma \ref{lem:entropyevol}, stating that
\begin{equation*}
    H(f_t \vert g_t) \leq H(f_0 \vert g_0) + \int_0^t (-\mathcal{G}_s+\mathcal{B}_s) ds.
\end{equation*}
Combining Lemmas \ref{lem:good} and \ref{lem:bad} to estimate $-\mathcal{G}_s+\mathcal{B}_s$, the relative Fisher informations cancel out and only the relative entropy remains:
\begin{align*}
    - \mathcal{G}_s+ \mathcal{B}_s
    &\leq C\left[(1+K_{g1}^2(s))(M_{f}+M_{g})+c_0^{-1}(K_{g1}^2(s)+K_{g3}^2(s)) (M^2_f+M^2_{g})\right]H(f_s\vert g_s)\\
    & \leq \mathcal{C}\left(K_{g1}^2(s)+K_{g3}^2(s)\right)H(f_s\vert g_s).
\end{align*}
We thus have
\begin{align*}
    H(f_t \vert g_t) &\leq H(f_0 \vert g_0) + \mathcal{C}\int_0^t \left(K_{g1}^2(s)+K_{g3}^2(s)\right)H(f_s\vert g_s) ds,
\end{align*}
and we want to apply Gronwall's lemma to conclude. We need to check that the integrand in the right hand side is finite. This is the case because $H(f_s)$ is bounded on $[0,T]$ by $H(f_0)$ and $t\mapsto \int f_t\ln g_t dv$ is absolutely continuous from the expression \eqref{eq:testlng}, so $H(f_s\vert g_s)$ is bounded and the other term is integrable. Hence Gronwall's lemma applies and gives the result.
\end{proof}

\section{Logarithmic bounds}
\label{sec:logarithmicbounds}
In this section we prove Theorem \ref{thm:gbounds} by showing that the strong solution to the Landau equation satisfies all the bounds necessary for Theorem \ref{thm:rentropy} to apply. In fact, we show something stronger: the bounds hold with constant-in-time $K_{gi}$ (in the notations of Section \ref{sec:weakstrong}). We first build the solution $g$ using Guillen-Silvestre's existence result \cite{GuillenSilvestre2023}. We then rely on the parabolic form of the Landau equation and use a maximum principle to propagate the Maxwellian bounds from above and below. Finally, still relying on the parabolic form, the Schauder estimates allow us to prove the logarithmic bounds on $g$ and its derivatives. Recall that we concentrate on the very soft potentials case: we fix $\gamma\in[-3,-2]$ in all this section.

\subsection{Building the solution}

We begin by building the strong solution $g$ we will work with.

\begin{prop}
\label{prop:solg}
    Let $g_0$ be as in Theorem \ref{thm:gbounds}. There exists a classical solution $g=(g_t)_{t\in[0,+\infty)}$ to the Landau equation \eqref{eq:landau}. For any $T>0$ and $\alpha\in(0,1)$, $g$ is bounded on $[0,T]\times \mathbb{R}^3$, and Hölder continuous: $[g]_{\alpha,[0,T]\times \mathbb{R}^3}<+\infty$.
\end{prop}

\begin{proof}
We want to apply \cite[Theorem 1.2]{GuillenSilvestre2023}, which states that there exists a classical solution to the Landau equation with initial data $g_0$. Among the hypothesis of this result, it only remains to show that the initial Fisher information
$$I(g_0):=\int \vert  \nabla \ln g_0 \vert^2 g_0dv$$
is finite. Using interpolation between $C^0$ and $C^{2+\delta}$ (recalled in Lemma \ref{lem:interpolation1}), it is not hard to see that, for some $k$,
$$  \vert \nabla \ln g_0 \vert  \leq C \jap{v}^k,$$
because of \eqref{eq:smoothcond} and $\vert \ln g_0 \vert \leq C \jap{v}^2$ (provided by the Maxwellian bounds).
Since $g_0$ has moments of any order, this shows the finiteness of $I(g_0)$. Hence by \cite[Theorem 1.2]{GuillenSilvestre2023}, there exists $g=(g_t)_{t\in[0,+\infty)}$ a classical solution of the Landau equation \eqref{eq:landau} with initial data $g_0$.

The solution in \cite{GuillenSilvestre2023} is built using a short time existence result (\cite[Theorem 2.1]{GuillenSilvestre2023}, adapted from \cite{HendersonSnelsonTarfulea2019}), which shows that it is moreover bounded on $[0,T]\times \mathbb{R}^3$ for any $T>0$ (in fact, the quantity $\Vert\jap{v}^k g_t \Vert_{ L^\infty}$ remains bounded, for some $k>0$, as long as the solution exists).

Since for any $\alpha\in(0,1)$, $[g_0]_{\alpha,\mathbb{R}^3}<+\infty$ (once again by interpolation), and for any $k$, $\jap{v}^k g_0 \in L^\infty$ thanks to the Maxwellian bound from above, \cite[Theorem 1.4]{HendersonSnelsonTarfulea2019} yields the Hölder continuity in both time and velocity. In particular $[g]_{\alpha,[0,T]\times \mathbb{R}^3}<+\infty$ (recall that this seminorm measures continuity in velocity only, uniformly in time).
\end{proof}

For the remainder of this section, $g$ will be the fixed solution provided by Proposition~\ref{prop:solg} and we will work on a fixed finite time interval $[0,T]$. We recall that the energy of $g$ is conserved and its entropy is decreasing over time.

\subsection{Ellipticity and propagation of Maxwellian bounds}
\label{ssec:maximum principle}

The Landau equation \eqref{eq:landau} can be written as a parabolic equation in nondivergence form:
\begin{equation}
    \label{eq:landaupara}
    \partial_t g = \bar{a} : \nabla^2g  +\bar{c}g.
\end{equation}
Written as is, it is a linear equation, but whose coefficients depend on $g$:
\begin{equation}
\label{eq:bara}
    \bar{a}(t,v) := \int_{\mathbb{R}^3} a(v-w)g_t(w)dw,
\end{equation}
\begin{equation}
\label{eq:barc}
\bar{c}(t,v) := - \sum_{i,j} \partial_{ij}^2\bar{a}_{ij}(t,v)= \begin{cases}
    2(\gamma + 3)\int_{\mathbb{R}^3} \vert v \vert^\gamma g_t(w)dw & \text{ if }\gamma \in (-3,-2] \\
    8\pi g_t(v) & \text{ if }\gamma=-3.
\end{cases}
\end{equation}
From now on, we will entirely rely on the formulation \eqref{eq:landaupara} rather than \eqref{eq:landau}. The goal of this section is to propagate the Maxwellian bounds on $g_0$ using the parabolic maximum principle. To make use of the parabolicity, we first need some classical ellipticity estimates on the symmetric matrix $\bar{a}$, that will also be useful in the next section:
\begin{lemma}[Ellipticity estimates, lower bound]
\label{lem:ellipticlower}
The symmetric matrix
$\bar{a}(t,v)$
satisfies for any time $t\in[0,T]$ and any unitary vector $\xi \in \mathbb{R}^3$
$$\langle \bar{a}(v) \xi , \xi\rangle \geq \lambda \jap{v}^{\gamma},$$
%and if in addition $\langle \xi , v \rangle= 0$,
%$$\langle \bar{a}(v) \xi , \xi\rangle \geq \lambda \jap{v}^{\gamma+2}.$$
where $\lambda>0 $ depends only on (upper bounds on) the energy and entropy of $g_0$.
\end{lemma}
This is exactly the version without cut-off of Lemma \ref{lem:ellipticitycutoff}. A proof can be found in \cite{Silvestre2016} for instance. We also have an upper bound:
\begin{lemma}[Ellipticity estimates, upper bound]
\label{lem:elliptichigher}
The symmetric matrix
$\bar{a}(t,v)$
satisfies for any time $t\in[0,T]$  and any unitary vector $\xi \in \mathbb{R}^3$
$$\langle \bar{a}(v) \xi , \xi\rangle \leq \Lambda \jap{v}^{\gamma+2},$$
and if $\xi$ is parallel to $v$,
$$\langle \bar{a}(v) \xi , \xi\rangle \leq \Lambda \jap{v}^{\gamma}.$$
Here $\Lambda>0 $ depends on (upper bounds on) $\Vert g \Vert_{L^\infty([0,T]\times\mathbb{R}^3)}$ and on the $p$-th moment of $g_0$, with $p=\frac{-3\gamma}{\gamma +5}-\gamma$.
\end{lemma}
\begin{remark}
    The value of $p$ ranges between $2$ for $\gamma=-2$ to $\frac{9}{2}$ for $\gamma=-3$. These estimates also hold for $\gamma>-2$ but their proof differs and does not require boundedness or additional moments.
\end{remark}
\begin{proof}
The proof is adapted from \cite[Lemma 2.1]{CameronSilvestreSnelson2018a} (in which the case $\gamma\in[-2,0)$ is treated) and from \cite[Lemma A3]{HendersonSnelson2019}. We fix a time $t\in[0,T]$ and drop the time dependence. We begin with the first point. Recalling that $a(v)=\vert v \vert^{\gamma+2}\Pi(v)$
\begin{align*}
    \langle \bar{a}(v)\xi , \xi \rangle &= \int_{\mathbb{R}^3}  \vert v-w \vert^{\gamma+2} \langle \Pi(v-w)\xi , \xi \rangle g(w)dw\\
    &\leq \int_{\mathbb{R}^3}  \vert v-w \vert^{\gamma+2} g(w)dw.
\end{align*}
Let $0<r \leq 1$ to be fixed later, we can cut this integral in three pieces (with the second one being possibly empty):
\begin{align*}
    \langle \bar{a}(v)\xi , \xi \rangle &\leq \int_{B(v,r)}  \vert v-w \vert^{\gamma+2} g(w)dw+\int_{B(v,\vert v \vert/2) \setminus B(v,r)}  \vert v-w \vert^{\gamma+2} g(w)dw\\
    &+\int_{\mathbb{R}^3 \setminus B(v,\max(r,\vert v \vert/2))}  \vert v-w \vert^{\gamma+2} g(w)dw.
\end{align*}
On the first one, we use the boundedness of $g$ to get:
$$\int_{B(v,r)}  \vert v-w \vert^{\gamma+2} g(w)dw \leq  \Vert g \Vert_{L^\infty(\mathbb{R}^3)} \int_{B(0,r)}  \vert w \vert^{\gamma+2} dw = C\Vert g \Vert_{L^\infty(\mathbb{R}^3)}r^{\gamma+5}.$$

The second one exists only when $\vert v \vert > 2r$. For $w\in B(v,\vert v \vert/2) \setminus B(v,r)$, we have $\vert v-w\vert \geq r$, as well as $\vert w \vert \geq \vert v \vert/2$. The latter implies that $\frac{1}{2} \leq \frac{\jap{w}}{\jap{v}}$, so that for any $0\leq k\leq  p$,
\begin{align*}
    \int_{B(v,\vert v \vert/2) \setminus B(v,r)}  \vert v-w \vert^{\gamma+2} g(w)dw &\leq 2^k r^{\gamma +2} \int_{B(v,\vert v \vert/2) \setminus B(v,r)}  \frac{\jap{w}^k}{\jap{v}^k}g(w)dw \\&\leq  \frac{2^k r^{\gamma +2}}{\jap{v}^k} C_g,
\end{align*}
with $C_g$ being the $p$-th moment of $g$.

On the last piece, $\vert v-w\vert \geq \max(r,\vert v \vert/2)$. We can find a small $c>0$ (independent of $r$) such that if $\vert v \vert /2 \geq 1$, then $\vert v\vert/2\geq c \langle v \rangle$, and else $r \geq r c\jap{v}$. Hence
$$\int_{\mathbb{R}^3 \setminus B(v,\max(r,\vert v \vert/2))}  \vert v-w \vert^{\gamma+2} g(w)dw \leq \mathcal{A}:=\begin{cases}
    c^{\gamma+2}  \jap{v}^{\gamma+2} &\text{ if }\vert v \vert /2 \geq 1\\
    c^{\gamma+2} r^{\gamma+2} \jap{v}^{\gamma+2} &\text{ else. }
\end{cases}$$
To summarize, 
\begin{align*}
    \langle \bar{a}(v)\xi , \xi \rangle &\leq C_g\left[ r^{\gamma+5} +\frac{4^k r^{\gamma +2}}{\jap{v}^k}\right] + \mathcal{A}
\end{align*}
We now choose $r=\jap{v}^\frac{\gamma+2}{\gamma +5}$, and $k=\frac{(\gamma+2)^2}{\gamma +5}-(\gamma +2) \in [0,\frac{3}{2}]$, such that
$$\frac{ r^{\gamma +2}}{\jap{v}^k} = \jap{v}^{\gamma+2}.$$
Then the first two terms are controlled by $\jap{v}^{\gamma+2}$. But $\mathcal{A}$ as well, because $r^{\gamma+2}$ is a power of $\jap{v}$, hence remains bounded for $\vert v \vert /2 \leq 1$.
\\

We now prove the second point. We can suppose $\xi = \frac{v}{\vert v \vert}$. Since we only want a better decay for large $v$, we can also suppose that $\vert v \vert \geq 2$. We proceed more carefully than before and write
\begin{align*}
    \langle \bar{a}(v)\xi , \xi \rangle &= \int_{\mathbb{R}^3}  \vert v-w \vert^{\gamma+2} \langle \pi(v-w)\xi , \xi \rangle g(w)dw\\
    &= \int_{\mathbb{R}^3}   \vert v-w \vert^{\gamma+2} \left(1 - \frac{\langle v-w , \xi \rangle^2}{\vert v-w\vert^2}\right) g(w)dw\\
    &=\int_{\mathbb{R}^3}   \vert v-w \vert^{\gamma} \left(\vert v-w\vert^2 - \langle v-w , \xi \rangle^2\right) g(w)dw\\
    &=\int_{\mathbb{R}^3}   \vert v-w \vert^{\gamma} \left(\vert v-w\vert^2 -(\vert v \vert- \langle w , \xi \rangle)^2\right) g(w)dw\\
\end{align*}
using that $\langle v, \xi \rangle = \vert v\vert$ in the last line. Expanding the square and simplifying,
\begin{align*}
    \langle \bar{a}(v)\xi , \xi \rangle &=\int_{\mathbb{R}^3}   \vert v-w \vert^{\gamma} \left(\vert w\vert^2 -\langle w , \xi \rangle^2\right) g(w)dw\\
    &=\int_{\mathbb{R}^3}   \vert v-w \vert^{\gamma} \vert w\vert^2 \sin^2(\theta) g(w)dw,
\end{align*}
with $\theta$ the angle between $v$ and $w$. We let $0<r\leq 1$ and cut the integral in three again. The first piece is the integral over $B(v,r)$. By hypothesis, $r \leq \vert v \vert / 2$, and for $w\in B(v, \vert v \vert /2)$, $\vert w \vert \vert \sin (\theta) \vert \leq \vert v-w \vert$. Hence
$$\int_{B(v,r)}   \vert v-w \vert^{\gamma} \vert w\vert^2 \sin^2(\theta) g(w)dw \leq\int_{B(v,r)}   \vert v-w \vert^{\gamma+2} g(w)dw \leq C\Vert g \Vert_{L^\infty(\mathbb{R}^3)}r^{\gamma+5}.$$
For $w\in B(v,\vert v \vert/2) \setminus B(v,r)$, $\vert w \vert \vert \sin (\theta) \vert \leq \vert v-w \vert$ still holds, and as in the first point we have $\vert v-w\vert \geq r$, as well as $\frac{1}{4} \leq \frac{\jap{w}}{\jap{v}}$, so that for any $0<k\leq p$,
$$\int_{B(v,\vert v \vert/2) \setminus B(v,r)}   \vert v-w \vert^{\gamma} \vert w\vert^2 \sin^2(\theta) g(w)dw \leq \frac{4^{k} r^{\gamma+2} }{ \jap{v}^k}\int_{\mathbb{R}^3}   \jap{w}^{k} g(w)dw \leq C_g\frac{r^{\gamma+2} }{ \jap{v}^{k}}.$$
Finally for $w\in \mathbb{R}^3 \setminus B(v,\vert v \vert/2)$, we drop the sine and simply write
$$\int_{\mathbb{R}^3 \setminus B(v,\vert v \vert/2)}   \vert v-w \vert^{\gamma} \vert w\vert^2  g(w)dw \leq \frac{\vert v\vert^{\gamma}}{2^\gamma} \int_{\mathbb{R}^3} \vert w \vert^2 g(w)dw \leq C_g\jap{v}^{\gamma}.$$
Choosing $r=\jap{v}^\frac{\gamma}{\gamma+5}$ and $k=\frac{\gamma(\gamma+2)}{\gamma +5}-\gamma = p$ concludes.

We have used the finiteness of the $p$-th moment of $g_t$, but by the propagation of moments (Lemma \ref{lem:momentbounds}), it can be bounded by the $p$-th moment of $g_0$ (which is finite because $g_0$ is bounded by a Maxwellian).
\end{proof}

We now bound the zero order term $\bar{c}$. The bound is blunt in the case $\gamma>-3$ but since we are mainly interested in the Coulomb case we do not wish to sharpen it.
\begin{lemma}
\label{lem:boundc}
The coefficient
$\bar{c}(t,v)$
is non-negative and bounded by a constant $\Vert \bar{c} \Vert_\infty$ depending on $\Vert g \Vert_{L^\infty([0,T]\times\mathbb{R}^3)}$.
\end{lemma}
\begin{proof}
    Recall that the solution $g$ given by Proposition \ref{prop:solg} is bounded on $[0,T]\times\mathbb{R}^3$. There is nothing to prove if $\gamma=-3$. Else, $\vert v \vert^\gamma$ is integrable around $0$, and for any time
    $$0 \leq \int_{\mathbb{R}^3} \vert v \vert^\gamma g_t(w)dw \leq \int_{\mathbb{R}^3} (\vert v \vert^\gamma \mathbf{1}_{\vert v \vert \leq1} + 1)g_t(w)dw \leq C(\Vert g_t \Vert_{L^\infty(\mathbb{R}^3)}+1),$$
    which concludes.
\end{proof}
We emphasize that the ellipticity constants $\lambda$ and $\Lambda$ and the bound on $\bar{c}$ are uniform in time. These estimates allow us to propagate the initial Maxwellian bounds. The proof relies on the parabolic maximum principle, the exact form we use being stated and proved in Appendix \ref{app:proofmaxprinc} for the convenience of the reader.
\begin{prop}[Propagation of Maxwellian bounds]
\label{prop:propbounds}
Recall that $g_0$ is bounded above and below by a Maxwellian:
$$ke^{-\mu\vert v \vert^2}\leq g_0(v) \leq Ke^{-\mu\vert v \vert^2},$$
for some $k,K,\mu>0$. This bound propagates: there exists $k',K'>0$ such that for all times $t\in[0,T]$
$$k'e^{-\mu\vert v \vert^2} \leq g_t(v) \leq K'e^{-\mu\vert v \vert^2}.$$
\end{prop}
\begin{proof}
Let $\psi = \exp(\omega t - \mu \vert v \vert^2)$. We want to make $\psi$ either a super- or subsolution of the parabolic equation \eqref{eq:landaupara} and then apply the maximum principle. A direct computation yields:
$$\nabla^2 \psi = (- 2\mu \Id + 4 \mu^2 v \otimes v) \psi,$$
so that dropping all non-negative terms and then using Lemma \ref{lem:elliptichigher},
\begin{align*}
    \bar{a} :\nabla^2 \psi + \bar{c} \psi &\geq \bar{a} :(-2\mu \Id)\psi\\
    &\geq -6\Lambda \mu \jap{v}^{\gamma+2} \psi\\
    &\geq -6\Lambda \mu\ \psi
\end{align*}
since $\gamma+2\leq 0$. Choosing $\omega=-6\Lambda \mu$, we make $\psi$ a subsolution to
$$\partial_t\psi \leq \bar{a} :\nabla^2 \psi + \bar{c} \psi.$$
Applying the parabolic maximum principle (Proposition~\ref{thm:maxprinciple}) to the subsolution $k\psi - g$ which is non-positive at $t=0$, we get that it remains non-positive, so that for all times $t\in[0,T]$
$$g_t(v) \geq k e^{\omega T} e^{- \mu \vert v \vert^2}.$$
Similarly, dropping the non-positive term and using the improvement in Lemma \ref{lem:elliptichigher} in directions parallel to $v$, as well as Lemma \ref{lem:boundc},
\begin{align*}
    \bar{a} :\nabla^2 \psi + \bar{c} \psi &\leq \left(4\mu^2 \bar{a} : v\otimes v + \Vert \bar{c} \Vert_\infty\right) \psi\\
    &\leq\left( 4\mu^2 \Lambda \jap{v}^\gamma \vert v \vert^2+ \Vert \bar{c} \Vert_\infty\right) \psi\\
    &\leq \left(4\mu^2 \Lambda+ \Vert \bar{c} \Vert_\infty\right) \psi.
\end{align*}
Choosing this time $\omega=4\mu^2 \Lambda+ \Vert \bar{c} \Vert_\infty$, we make $\psi$ a supersolution and apply Proposition~\ref{thm:maxprinciple} to $g - K \psi$ to conclude.
\end{proof}
\begin{remark}
    In the previous proposition it is crucial that $\gamma\leq-2$, otherwise we would only be able to propagate a \textit{deteriorating} Maxwellian, i.e. bound $g$ from below and above as
    $$e^{-(\mu+\varepsilon t)\vert v \vert^2}\lesssim g_t \lesssim e^{-(\mu-\varepsilon t)\vert v \vert^2}$$
    for any small $\varepsilon>0$. 
\end{remark}

\subsection{Using Schauder estimates to derive logarithmic bounds}

We now turn to the estimation of the derivatives of $\ln g$. We will obtain them by estimating the derivatives of $g$ in terms of $g$ itself, which after dividing everything by $g$ will yield estimates on $\ln g$. Since $g$ satisfies the parabolic equation \eqref{eq:landaupara}, we can control its derivatives using parabolic Schauder estimates.
To apply these Schauder estimates, we need Hölder continuity of 
 the coefficients $\bar{a}$ and $\bar{c}$. This is done in the following two lemmas, the first focusing on the coefficient $\bar{a}$.
\begin{lemma}[Hölder continuity of $\bar{a}$]
\label{lem:Höldera}
    The matrix-valued function $\bar{a}(t,v)$ is $\alpha$-Hölder continuous in $v$ on $\mathbb{R}^3$, uniformly in $t$, with
    $\alpha = \frac{\gamma+5}{4}$. Moreover,
    $$[\bar{a}]_{\alpha,[0,T]\times\mathbb{R}^3} \leq C(1+\Vert g\Vert_{L^\infty([0,T]\times\mathbb{R}^3)}),$$
    for $C>0$ depending only on $\alpha$.
\end{lemma}
\begin{proof}
    Let us consider a fixed $t\in[0,T]$, and simply write $g$ for $g_t$. For $v,v'\in \mathbb{R}^3$, we write:
    \begin{align*}
        \Vert \bar{a}(v)-\bar{a}(v') \Vert \leq \int_{\mathbb{R}^3}&\left\Vert \vert v-w \vert^{\gamma+2} \Pi(v-w)-\vert v'-w \vert^{\gamma+2}\Pi(v'-w)\right\Vert g(w)dw\\
        \leq \int_{\mathbb{R}^3}&\left\vert \vert v-w \vert^{\gamma+2} -\vert v'-w \vert^{\gamma+2}\right\vert g(w)dw \\
        &+ \int_{\mathbb{R}^3}\vert v'-w \vert^{\gamma+2}\left\Vert  \Pi(v-w)-\Pi(v'-w)\right\Vert g(w)dw \\
        = \int_{\mathbb{R}^3}&\left\vert \vert v-v'+w\vert^{\gamma+2} -\vert w \vert^{\gamma+2}\right\vert g(v'-w)dw \\ &+\int_{\mathbb{R}^3}\vert w \vert^{\gamma+2}\left\Vert  \Pi(v-v'+w)-\Pi(w)\right\Vert g(v'-w)dw.
    \end{align*}
using that $\Vert \Pi \Vert \leq 1$.

Let us begin with the first term, for which we can suppose $\gamma< -2$. Fix $r>0$ to be determined later. Since
for any $z\in \mathbb{R}^3\setminus B(0,r)$,  $\vert \nabla\vert z \vert^{\gamma+2} \vert \approx \vert z \vert^{\gamma+1}\leq  r^{\gamma+1}$, using some path between $z$ and $z'$ of length $\leq C \vert z - z'\vert$ staying outside of $B(0,r)$ (which always exist, and we emphasize that $C$ does not depend on $r$), we get the following Lipschitz behavior: for all $z,z'\in \mathbb{R}^3\setminus B(0,r)$,
$$ \left\vert \vert z\vert^{\gamma+2} - \vert z'\vert^{\gamma+2} \right\vert \leq Cr^{\gamma+1}\vert z - z'\vert. $$
Let $U=B(0,r)\cup B(v'-v,r)$. On $\mathbb{R}^3 \setminus U$ we can plug the above bound,
$$\int_{\mathbb{R}^3 \setminus U}\left\vert \vert v-v'+w\vert^{\gamma+2} -\vert w \vert^{\gamma+2}\right\vert g(v'-w)dw \leq Cr^{\gamma+1}\left\vert v-v'\right\vert,  $$
using that $g$ has unit mass.

For the remaining integral over $U$, we simply write
\begin{align*}
    \int_{U}\left\vert \vert v-v'+w\vert^{\gamma+2} -\vert w \vert^{\gamma+2}\right\vert g(v'-w)dw &\leq \Vert g \Vert_{L^\infty} \int_{ U} \left(\vert v-v'+w\vert^{\gamma+2} +\vert w \vert^{\gamma+2}\right)dw.
\end{align*}
Because $U$ is the union of two balls, the right hand side is bounded by a sum of integrals of $w \mapsto \vert w \vert^{\gamma+2}$ over balls of radius $r$ and varying centers. But the largest possible value for such an integral is when the center is $0$: \textit{i.e.} for any center $v_0$, 
$$ \int_{B(v_0,r)} \vert w \vert^{\gamma+2} dw \leq \int_{B(0,r)} \vert w \vert^{\gamma+2} dw = Cr^{\gamma+5}.$$
Hence combining both inequalities
$$\int_{\mathbb{R}^3}\left\vert \vert v-v'+w\vert^{\gamma+2} -\vert w \vert^{\gamma+2}\right\vert g(v'-w)dw \leq C(1+\Vert g \Vert_{L^\infty})\left[r^{\gamma+1}\left\vert v-v'\right\vert+ r^{\gamma+5}\right].$$
Recall that $\alpha = \frac{\gamma+5}{4}$, and choose $r=\vert v-v'\vert^{\frac{\alpha-1}{\gamma +1}}$. It yields $r^{\gamma+5}=\vert v-v'\vert^\alpha$ because
$$\frac{(\alpha-1)(\gamma+5)}{\gamma +1}=\frac{(\frac{\gamma+1}{4})(\gamma+5)}{\gamma +1}=\alpha,$$
and this term is $\alpha$-Hölder continuous.

For the second term, we argue similarly. Using that the function $\Pi(z)$ depends only on $\frac{z}{\vert z \vert}$, it is easy to see that it is Lipschitz on $\mathbb{R}^3 \setminus B(0,1)$. Then, for any $r>0$, and $z,z'\in \mathbb{R}^3 \setminus B(0,r)$,
$$\Vert \Pi(z) - \Pi(z') \Vert = \left\Vert \Pi\left(\frac{z}{r}\right) - \Pi\left(\frac{z'}{r}\right) \right\Vert \leq \frac{C}{r} \vert z - z' \vert.$$
Using the same $U=B(0,r)\cup B(v'-v,r)$ to cut the integral in two as before, and also using the bound $\vert w\vert^{\gamma+2}\leq r^{\gamma+2}$ on $\mathbb{R}^3 \setminus U$, we get
$$\int_{\mathbb{R}^3}\vert w \vert^{\gamma+2}\left\Vert  \Pi(v-v'+w)-\Pi(w)\right\Vert g(v'-w)dw \leq C(1+\Vert g \Vert_{L^\infty})\left[r^{\gamma+1}\left\vert v-v'\right\vert+ r^{\gamma+5}\right],$$
which is the same bound as for the second term, so this concludes. The uniformity in time follows from the uniform boundedness of $g$.
\end{proof}
We now turn to the Hölder continuity of $\bar{c}$. In the Coulomb case, $\bar{c}$ is a multiple of $g$ so it was already shown back in Proposition \ref{prop:solg} that it is $\alpha$-Hölder continuous for any $\alpha\in(0,1)$. For the non-Coulomb case, we have:
\begin{lemma}[Hölder continuity of $\bar{c}$ (non-Coulomb case)]
\label{lem:Hölderc}
    For $\gamma \in ]-3,-2]$, the function $\bar{c}(t,v)$ is $\alpha$-Hölder continuous in $v$, uniformly in $t$, with
    $\alpha = \frac{\gamma+3}{4}$. Moreover, 
    $$[\bar{c}]_{\alpha,[0,T]\times\mathbb{R}^3} \leq C(1+\Vert g\Vert_{L^\infty([0,T]\times\mathbb{R}^3)}),$$
    for $C>0$ depending only on $\alpha$.
\end{lemma}
\begin{proof}
    We have
    $$\vert\bar{c}(t,v)-\bar{c}(t,v')\vert\leq \int_{\mathbb{R}^3}\left\vert \vert v-w\vert^{\gamma} -\vert v'-w \vert^{\gamma}\right\vert g_t(w)dw.$$
    The proof is thus exactly the same as the bound on the first term in the proof of Lemma~\ref{lem:Höldera} above, but $\gamma$ replaces $\gamma+2$, hence the new $\alpha$ exponent.
\end{proof}

With the Hölder continuity of the coefficients $\bar{a}$ and $\bar{c}$ in hand, we can apply Schauder estimates to the Landau equation \eqref{eq:landaupara}. Although parabolic Schauder estimates are a classical result, we are not in the most standard setting: our equation is not uniformly parabolic in $v$ (since the ellipticity degenerates as $v\rightarrow +\infty$), and we work with Hölder continuity in velocity only rather than time and velocity. We also wish to \textit{propagate} the Hölder continuity from the initial time rather than \textit{generate} it, \textit{i.e.} we will work in a \textit{boundary} case in time but an \textit{interior} case in velocity. For these reasons, and because we want to carefully keep track of the constants (especially the powers of $\jap{v}$), we recall a proof of the exact Schauder estimate we use in Appendix \ref{app:schauder}. An early reference on Schauder estimates in space only is \cite{Brandt1969}.

To state this Schauder estimate, we need some further notation. For the remainder of this section, we let $$\Sigma_R(v_0):=v_0+(-R,R)^3$$ be the cube in velocity centered on $v_0$ of side $2R$.
%It is in fact easier to work with $\vert \cdot \vert$ being the $\ell^\infty$ norm in $\mathbb{R}^3$, whose ball is the above cube, but every result can of course be transposed to the Euclidean norm by adjusting the constants.
We define
$$Q_R(v_0):=[0,T]\times \Sigma_R(v_0),$$
the cylinder in time and velocity.
We now need to define weighted norms. Let $v_0\in \mathbb{R}^3$, $v\in \Sigma_R(v_0)$, we define the distance to the boundary:
$$d_v:=\min\{\vert v-w\vert, w\in \partial \Sigma_R(v_0)\},$$
and for two points $d_{vw}=\min (d_v,d_w)$. The weighted Hölder norm is
$$[\varphi]^*_{\alpha,Q_{R}(v_0)}:=\sup_{\substack{(t,v), (t,w) \in Q_R(v_0)\\
v\neq w}} d_{vw}^\alpha\frac{\vert \varphi(t,v)-\varphi(t,w)\vert}{\vert v-w \vert^\alpha},$$
and the $2+\alpha$ version is
$$[\varphi]^*_{2+\alpha,Q_{R}(v_0)}:=\sup_{i,j}\sup_{\substack{(t,v), (t,w) \in Q_R(v_0)\\
v\neq w}} d_{vw}^{2+\alpha}\frac{\vert \partial_{ij}\varphi(t,v)-\partial_{ij}\varphi(t,w)\vert}{\vert v-w \vert^\alpha}.$$
By boundedness of $Q_R(v_0)$ and hence of $d_\cdot$, these semi-norms are weaker than the unweighted ones, but "equivalent far from the boundary". We also define $[\bar{c}]_\alpha^\bigstar$ as $[\bar{c}]^*_{\alpha,Q_{R}(v_0)}$ but with a power $2+\alpha$ instead of $\alpha$ in the weight. In those weighted norms, the Schauder estimate for the Landau equation writes:

\begin{prop}[Schauder estimate for the parabolic Landau equation]
\label{thm:schauderlandau}
    Let $v_0\in \mathbb{R}^3$ and $R\leq\frac{1}{2}\jap{v_0}$. Consider $u:Q_R(v_0)\rightarrow \mathbb{R}$, $C^1$ in time and $C^2$ in velocity, a solution of
    \begin{equation*}
        \partial_t u =\bar{a}:\nabla^2 u + \bar{c}u
    \end{equation*}
    with $\bar{a}$, $\bar{c}$ defined as in \eqref{eq:bara}, \eqref{eq:barc} (with $g$ from Proposition \ref{prop:solg}).

    Let $\alpha\in(0,1)$. If $u$ is $C^{2+\alpha}$ on $Q_{R}(v_0)$, it holds that:
    \begin{equation}
    \label{eq:schauderlandau}
        [u]^*_{2+\alpha,Q_{R}(v_0)} \leq C_{\alpha,\lambda,\Lambda} \left([u_0]^*_{2+\alpha,\Sigma_R(v_0)} + \mathcal{V}(v_0) \Vert u \Vert_{L^\infty(Q_R(v_0))}\right)
    \end{equation}    
    with \begin{align*}
        \mathcal{V}(v_0)=(1+[\bar{a}]^*_{\alpha,Q_{R}(v_0)})^{1+\frac{2}{\alpha}}&\jap{v_0}^{(2-\gamma)(1+\frac{2}{\alpha})}\\&+[\bar{c}]_\alpha^\bigstar\jap{v_0}^{-\gamma}+\left(1+\Vert d^2\bar{c}\Vert_{L^\infty(Q_R(v_0))}\right)^{1+\frac{\alpha}{2}}\jap{v_0}^{-\gamma(1+\frac{\alpha}{2})},
    \end{align*}
     where the constants $\lambda, \Lambda$ are from Lemmas \ref{lem:ellipticlower} 
    and \ref{lem:elliptichigher}.
\end{prop}
\begin{remark}
    At this point, we only know the Hölder continuity of $\bar{a}$ and  $\bar{c}$ for $0<\alpha \leq \frac{\gamma+3}{4}$ if $\gamma \in ]-3,-2]$ and $0<\alpha \leq \frac{1}{2}$ if $\gamma=-3$.  The Hölder norms of the coefficients $\bar{a}, \bar{c}$ above are possibly infinite, in which case the result is empty. We keep a general $\alpha$ because we will later improve on the Hölder continuity of $\bar{a}$ and  $\bar{c}$ by bootstrapping the results on $g$.
\end{remark}
In the exact expression of $\mathcal{V}(v_0)$, what matters most is the dependence on $\jap{v_0}$ because it will ultimately dictate the exponents that will appear in our estimates on $\ln g$ and its derivatives (which in turn dictate how many moments are required for the weak solution in the weak strong principle). We give the proof of this result in Appendix \ref{app:schauder} mainly in order to justify the expression of $\mathcal{V}(v_0)$, which is usually less precise in the literature. As it is usual with Schauder estimates, the proof consists in first proving a version for \textit{frozen coefficients} (that is to say, for constant $\bar{a}$ and $\bar{c}=0$), and then using the Hölder continuity of these coefficients to derive the general case.

\begin{remark}
    In \cite{HendersonSnelson2019, HendersonSnelsonTarfulea2019}, Schauder estimates for the \emph{inhomogeneous} Landau equation are obtained thanks to a change of variable that maps the parabolic Landau equation on (the equivalent of) the cylinder $Q_R(v_0)$ to a uniformly parabolic equation (with ellipticity constants uniform with respect to $v_0$), to which one can apply uniform Schauder estimates. A polynomial weight in $v_0$ then appears when doing the inverse change of variable. Since the works \cite{HendersonSnelson2019, HendersonSnelsonTarfulea2019} are not about tracking the exact polynomial weight, it is not clear if this change of variable method (which of course also applies to the simpler homogeneous setting) provides an improvement of the final result. However, it could be expected, because the change of variable takes into account the anisotropic improvement in the ellipticity estimates (see Lemma \ref{lem:elliptichigher}), while our method does not. We did not implement this method to keep this section a bit less technical.
\end{remark}

In the remaining of this section, we apply the Schauder estimate to $g$ itself to obtain the logarithmic bounds through some technical manipulations. We need to show that we can control the seminorm of the initial condition $[g_0]^*_{2+\alpha,\Sigma_R(v_0)}$ in the right hand side of the estimate~\eqref{eq:schauderlandau} by $\Vert g_0 \Vert_{L^\infty(\Sigma_R(v_0))}$ using the hypothesis on the initial condition \eqref{eq:smoothcond}. This is done is the following technical lemma, the proof of which is postponed to the end of the section. Recall that $\delta$ is such that $\ln g_0$ is locally $C^{2+\delta}$ and $\ell$ is the exponent in the polynomial growth of the $C^{2+\delta}$ seminorm (see \eqref{eq:smoothcond}).

\begin{lemma}
\label{lemma:reworksmoothcond}
    Let $\alpha\in(0,\delta]$. Let $v\in\mathbb{R}^3$ $r=\frac{1}{2}\jap{v}^{-\beta}$ with $\beta = \max(2+\frac{1}{2+\delta}(\ell-2), 1)$. It holds that
    \begin{equation}
        [g_0]^*_{2+\alpha,\Sigma_r(v)} \leq C \Vert g_0\Vert_{L^\infty( \Sigma_r(v))}\jap{v}^{-\beta \alpha+\sigma},
    \end{equation}
    with $\sigma=\max(2+\frac{\alpha}{2+\delta}(\ell-2),2-\alpha)$, for some $C>0$ depending on the constant $C_g$ in \eqref{eq:smoothcond}.
\end{lemma}

The values of $\beta$ and $\sigma$ are the results of the interpolations between norms done in the proof. Admitting this lemma, we are ready to give the proof of the following:

\begin{prop}[Logarithmic bounds]
\label{prop:logarithmic bounds}
    There exists a constant $K>0$ such that, for all $t\in[0,T]$:
    \begin{align*}
        \vert \nabla \ln g_t \vert &\leq K\jap{v}^\kappa\\
        \Vert \nabla^2 \ln g_t \Vert &\leq K\jap{v}^\zeta\\
        \vert \partial_t \ln g_t \vert &\leq K\jap{v}^\nu.
    \end{align*}
    The exponents are given by
    \begin{align*}
        \kappa &= 2 \beta+\frac{2-\gamma}{\delta},\\
        \zeta &= 2\kappa,\\
        \nu &=2\beta + \frac{4-2\gamma}{\delta},
    \end{align*}
    where $\beta = \max(2+\frac{1}{2+\delta}(\ell-2), 1)$ and $\delta$ is from the hypothesis of Theorem \ref{thm:gbounds}.
\end{prop}

\begin{proof}
    \textit{Strategy.} The proof is technical, mostly because we track the exact exponents $\kappa, \zeta, \nu$. We first describe the broad strategy. In Step 1, we apply the Schauder estimate on a narrow cylinder $Q_r(v)$ (imagine $v$ is large and $r$ is small) to obtain control of the derivatives
    $$ [ g ]_{2+\alpha, Q_r(v)} \lesssim \jap{v}^k  (\Vert g \Vert_{L^\infty(Q_r(v))} + [ g_0 ]_{2+\alpha, \Sigma_r(v)})$$
    for some $k,\alpha>0.$ Using Lemma \ref{lemma:reworksmoothcond}, we can control the right hand side by $$\jap{v}^m\Vert g \Vert_{L^\infty(Q_r(v))},$$ for some $m$. In Step 2, we show that the Maxwellian bounds allow us to write $$\Vert g \Vert_{L^\infty(Q_r(v))}\lesssim g_t(v).$$
    In Step 3, we interpolate to get
    $$\nabla^2 g (v),\nabla g (v) \lesssim \jap{v}^n  g_t(v),$$
    for some $n$ and this yields the logarithmic estimates by dividing by $g_t(v)$. It yields the result but with worst exponents $\kappa, \zeta, \nu$. In Step 4, we bootstrap the regularity obtained on $g$ to improve the exponents.
    \\
    
    \textit{Step 1.} We let $\alpha = \min(\frac{\gamma+3}{4},\delta)$ if $\gamma \in ]-3,-2]$ and $\alpha=\min(\frac{1}{2},\delta)$ if $\gamma=-3$. With this choice of $\alpha$, $\bar{a}$ and $\bar{c}$ are $\alpha$-Hölder continuous over the whole $\mathbb{R}^3$ thanks to Lemmas \ref{lem:Höldera},~\ref{lem:Hölderc} (and the boundedness of  $\bar{a}$ and $\bar{c}$, Lemmas \ref{lem:elliptichigher} and~\ref{lem:boundc}).

    Let $v\in\mathbb{R}^3$, and $r=\frac{1}{2}\jap{v}^{-\beta}$ (with $\beta = \max(2+\frac{1}{2+\delta}(\ell-2), 1)$ in view of applying Lemma \ref{lemma:reworksmoothcond}). Since $g$ solves the Landau equation,
    $$\partial_t g = \bar{a}:\nabla^2 g + \bar{c} g$$
    we can apply the Schauder estimate (Proposition \ref{thm:schauderlandau}) on the cube $Q_{r}(v)$. It yields
    \begin{equation*}
        [g]^*_{2+\alpha,Q_{r}(v)} \leq C\left[[g_0]^*_{2+\alpha,\Sigma_r(v)} + \mathcal{V}(v) \Vert g \Vert_{L^\infty(Q_{r}(v))}\right]
    \end{equation*}    
    with $$\mathcal{V}(v)=(1+[\bar{a}]^*_{\alpha,Q_{r}(v)})^{1+\frac{2}{\alpha}}\jap{v}^{(2-\gamma)(1+\frac{2}{\alpha})}+[\bar{c}]_\alpha^\bigstar\jap{v}^\gamma +\left(1+\Vert d^2\bar{c}\Vert_{L^\infty(Q_r(v))}\right)^{1+\frac{\alpha}{2}}\jap{v}^{-\gamma(1+\frac{\alpha}{2})}.$$ 
     Using that the weight in the Hölder norms are $\lesssim r \leq 1$, the norms can be bounded independently of $v$ and $r$, so we get:
    $$\mathcal{V}(v) \leq C \jap{v}^{(2-\gamma)(1+\frac{2}{\alpha})}.$$
    Plugging in Lemma \ref{lemma:reworksmoothcond} (and bounding $\Vert g_0\Vert_{L^\infty( \Sigma_r(v))}\leq \Vert g\Vert_{L^\infty( Q_r(v))}$), we get
    \begin{align}
    \label{eq:schauderstep1}
        [g]^*_{2+\alpha,Q_{r}(v_0)} &\leq C\left[ \jap{v}^{-\beta\alpha+\sigma}+\jap{v}^{(2-\gamma)(1+\frac{2}{\alpha})}\right]\Vert g\Vert_{L^\infty( Q_r(v))}\nonumber \\
        &\leq C\jap{v}^{(2-\gamma)(1+\frac{2}{\alpha})}\Vert g\Vert_{L^\infty( Q_r(v))},
    \end{align}
    because since $\sigma= \max(2+\frac{\alpha}{2+\delta}(\ell-2), 2-\alpha)$, we always have 
    $$-\beta \alpha + \sigma =2-2\alpha \leq 2 \leq  (2-\gamma)(1+\frac{2}{\alpha}).$$
    The right hand side is in terms of $\Vert g\Vert_{L^\infty( Q_r(v))}$, which is what we wanted.
    \\
    
    \textit{Step 2.} In this step, we show that $\Vert g\Vert_{L^\infty( Q_r(v))}$ can be controlled by the value $g$ takes at the center at any point $t$ in time. It holds that
    $$ - ( \vert v \vert -r)^2   = -\vert v \vert^2 +\vert v \vert \jap{v}^{-s}-\frac{1}{4}\jap{v}^{-2s}  \leq 1-\vert v \vert^2,$$
    so that
    $$\sup_{w\in Q_r(v)} e^{-\mu \vert w \vert^2} = e^{-\mu (\vert v \vert-r)^2} \leq e^{\mu-\mu \vert v \vert^2}.$$
    Using the propagated Maxwellian bounds (Propositon \ref{prop:propbounds})
    \begin{equation}
    \label{eq:step1}
        \Vert g\Vert_{L^\infty( Q_r(v))}\leq K'\sup_{w\in Q_r(v)} e^{-\mu \vert w \vert^2} \leq K'e^{\mu-\mu \vert v \vert^2} \leq \left(\frac{K'}{k'}e^{\mu} \right) g_t(v).
    \end{equation}
    This is the control we sought.
    \\

    \textit{Step 3.} In this step we want a control of
    $\nabla g_t (v)$ and $\nabla^2 g_t(v)$. What we control for the moment is the $C^{2+\alpha}$ seminorm of $g$, so we want to drop down by interpolation. We work on the weighted norms using Lemma \ref{lem:interpolation2}, with $\varepsilon= \jap{v}^{-\eta}$, for some $\eta\geq 0$. For any time $t$,
    \begin{align*}
        \Vert \nabla^2 g_t(v) \Vert&\leq r^{-2}\Vert d^2 \nabla^2 g\Vert_{L^\infty(Q_{r}(v))}\\
        &\leq C r^{-2}\left(\jap{v}^{-\alpha \eta}[g]^*_{2+\alpha,Q_{r}(v)} + \jap{v}^{2\eta}\Vert g \Vert_{L^\infty(Q_{r}(v))}\right)\\
        &\leq C\left(\jap{v}^{(2-\gamma)(1+\frac{2}{\alpha}) +2\beta-\alpha \eta}+\jap{v}^{2 \beta+2\eta}\right)\Vert g \Vert_{L^\infty(Q_{r}(v))},
    \end{align*}
    where the second line is obtained by interpolation and the third by plugging in the Schauder estimate from Step 1 \eqref{eq:schauderstep1}.
    Choosing $\eta=\frac{2-\gamma}{\alpha}$, we get:
    \begin{align*}
        \Vert \nabla^2 g_t(v) \Vert
        &\leq C\jap{v}^{2 \beta+\frac{4-2\gamma}{\alpha}}\Vert g \Vert_{L^\infty(Q_{r}(v))}\\
        &\leq C\jap{v}^{2 \beta+\frac{4-2\gamma}{\alpha}}g_t(v),
    \end{align*}
    using \eqref{eq:step1} from Step 2 for the last line.

    We proceed similarly for the gradient. The interpolation writes:
    \begin{align*}
        \vert \nabla g_t(v) \vert&\leq r^{-1}\Vert d \nabla g\Vert_{L^\infty(Q_{r}(v))}\\
        &\leq C r^{-1}\left(\jap{v}^{-(1+\alpha) \eta}[g]^*_{2+\alpha,Q_{r}(v)} + \jap{v}^{\eta}\Vert g \Vert_{L^\infty(Q_{r}(v))}\right)\\
        &\leq C\left(\jap{v}^{(2-\gamma)(1+\frac{2}{\alpha}) +2\beta-(1+\alpha)\eta}+\jap{v}^{2 \beta+\eta}\right)\Vert g \Vert_{L^\infty(Q_{r}(v))}.
    \end{align*}
    Choosing once again $\eta=\frac{2-\gamma}{\alpha}$,
    \begin{align}
    \label{eq:graddecay}
        \vert \nabla g_t(v) \vert
        &\leq C\jap{v}^{2 \beta+\frac{2-\gamma}{\alpha}}g_t(v).
    \end{align}
    From this we obtain the first logarithmic bound we wanted, the one for the gradient:
    $$ \vert \nabla \ln g_t(v) \vert \leq C\jap{v}^{2 \beta+\frac{2-\gamma}{\alpha}},$$
    and using
    $$\nabla^2 \ln g = \frac{\nabla^2 g}{g} -\nabla \ln g \otimes \nabla \ln g,$$
    we get the logarithmic Hessian bound:
    $$ \Vert \nabla^2 \ln g_t(v) \Vert \leq C\jap{v}^{4 \beta+\frac{4-2\gamma}{\alpha}}.$$
    
    Finally, using the Landau equation with the bound $\Vert \bar{a} \Vert \lesssim  \Lambda \jap{v}^{\gamma+2} \leq \Lambda$ and the bound $\bar{c}\leq \Vert \bar{c} \Vert_{L^\infty([0,T]\times\mathbb{R}^3)}$, we also get
    $$\vert \partial_t g(v)\vert \leq C\jap{v}^{2 \beta+\frac{4-2\gamma}{\alpha}} g_t(v).$$
    We have estimated every quantity we sought. However, the exponents we obtained feature $\alpha$ in their denominator, rather than $\delta$ as in the statement of the proposition (and $\alpha\leq \delta$ so this is indeed worse). The next step shows how to improve them.
    \\
    
    \textit{Step 4.}  We can do slightly better by using the estimates from Step 3 to show that $\bar{a}$ and $\bar{c}$ are in fact Lipschitz in $v$, uniformly in time. The estimate \eqref{eq:graddecay} together with the Maxwellian upper bound for $g$ show that
    $$\vert \nabla g_t(v) \vert
        \leq Ce^{-\frac{\mu }{2}\vert v \vert^2}.$$
    Plugging this back in the expression of $\bar{a}$, this means that for any $v,v'$ such that $\vert v-v'\vert \leq 1$,
    \begin{align*}
        \Vert \bar{a}(t,v) - \bar{a}(t,v')\Vert &\leq \int_{\mathbb{R}^3} \vert w \vert^{\gamma+2} \vert g_t(v-w)-g_t(v'-w)\vert dw\\
        &\leq\vert v-v'\vert\int_{\mathbb{R}^3} \vert w \vert^{\gamma+2}  \Vert \nabla g_t \Vert_{L^\infty([v-w,v'-w])} dw.
    \end{align*}
    But $z\in [v-w,v'-w]$ satisfies $\vert z \vert \geq  \vert v-w\vert - \vert v-v'\vert\geq \vert v-w\vert -1$, so that
    $$\Vert \nabla g_t \Vert_{L^\infty([v-w,v'-w])} \leq Ce^{-\frac{\mu}{2}(\vert v-w\vert -1)^2_+}$$
    and one easily sees that
    $$\sup_v \int_{\mathbb{R}^3} \vert w \vert^{\gamma+2}  e^{-\frac{\mu}{2}(\vert v-w\vert -1)^2_+} dw <+\infty.$$
    Together with the boundedness of $\bar{a}$, this shows that $\bar{a}$ is Lipschitz. The same proof works on $\bar{c}$ in the $\gamma>-3$ case, so in any case both coefficients are Lipschitz. In particular they are $\alpha$-Hölder continuous for any $\alpha\in(0,1)$. We can thus repeat Steps 1, 2 and 3, but choosing $\alpha=\delta$ at the beginning of Step 1. We end up with the bounds above but with $\delta $ instead of $ \alpha$, and $\beta$ unchanged.
\end{proof}

We can finally give the proof of Theorem \ref{thm:gbounds}.
\begin{proof}\textit{(Theorem \ref{thm:gbounds})}
    The solution $g$ given by Proposition \ref{prop:solg} satisfies all the conditions of Theorem \ref{thm:rentropy}. Indeed, it is qualitatively smooth on $[0,T]\times\mathbb{R}^3$ and the logarithmic bounds hold in $L^\infty$ in time by Proposition \ref{prop:logarithmic bounds}. For any $f_0$ with sufficient moments, the relative entropy estimate holds.
\end{proof}

All that remains to prove is the technical Lemma \ref{lemma:reworksmoothcond} on the initial condition:

\begin{proof} \textit{(Lemma \ref{lemma:reworksmoothcond})}

Recall that the goal is to  control $[ g_0 ]^*_{2+\alpha, \Sigma_r(v)}$ in terms of $\Vert g_0 \Vert_{L^\infty(\Sigma_r(v))}$ and some power of $\jap{v}$, using the hypothesis \eqref{eq:smoothcond}. We will in fact work with unweighted norms, and add the weight at the very end.

Relying on the formula
    \begin{equation}
    \label{eq:formulaforstep2}
        \frac{\nabla^2 g_0}{g_0} = \nabla^2 \ln g_0+\nabla \ln g_0 \otimes \nabla\ln g_0,
    \end{equation}
    we first see how each term in the right hand side grows for large velocities in $L^\infty$ and $C^\alpha$, using interpolations on a large cube $\Sigma_R(0)$ (Step 1). We then write $\nabla^2 g_0 = g_0 \frac{\nabla^2 g_0}{g_0}$ to conclude with some more interpolations on $\Sigma_r(v)$ (Step 2).
    
    \textit{Step 1.} Let $R\geq 1$. We use the interpolation Lemma \ref{lem:interpolation1} on $\Sigma_R(0)$ with $\varepsilon=\frac{1}{2}R^{-\eta}$. For the condition $\varepsilon<R$ required by this lemma to hold, $\eta$ must satisfy $\eta\geq -1$. The interpolation writes:
    \begin{align*}
        \Vert \nabla \ln g_0 \Vert_{L^\infty(\Sigma_R(0))} &\leq C \left( R^{ - \eta(1+\delta)}[\ln g_0 ]_{2+\delta,\Sigma_R(0)} + R^{\eta}\Vert \ln g_0 \Vert_{L^\infty(\Sigma_R(0))}\right)\\
        &\leq C \left( R^{\ell - \eta(1+\delta)} + R^{2+\eta}\right)
    \end{align*}
    where the second line was obtained using $$\Vert \ln g_0 \Vert_{L^\infty(\Sigma_R(0))} \lesssim R^2$$ (obtained from the Maxwellian bound) and $$[\ln g_0 ]_{2+\delta,\Sigma_R(0)} \lesssim R^\ell,$$ which is the hypothesis \eqref{eq:smoothcond} written with cubes instead of balls. Choosing $\eta =\frac{\ell-2}{2+\delta}$ if $\ell> -\delta$, and $\eta=-1$ otherwise, this gives
    \begin{equation}
    \label{eq:step2gradcontrol}
        \Vert \nabla \ln g_0 \Vert_{L^\infty(\Sigma_R(0))} \leq CR^{\beta},
    \end{equation}
    with $\beta = \max(2+\frac{1}{2+\delta}(\ell-2), 1)$ as announced.

    Similarly,
    $$[ \nabla \ln g_0 ]_{\alpha, \Sigma_R(0)} \leq C \left( R^{\ell - \eta(1+\delta-\alpha)} + R^{2+\eta(1+\alpha)}\right) $$
    so
    $$[ \nabla \ln g_0 ]_{\alpha, \Sigma_R(0)} \leq CR^{\theta},$$
    with $\theta = \max(2+\frac{1+\alpha}{2+\delta}(\ell -2), 1-\alpha)$ (obtained with the same choice of $\eta$).

    Finally
    $$[\nabla^2\ln g_0]_{\alpha, \Sigma_R(0)} \leq C R^\epsilon$$
    with $\epsilon=\max(2+\frac{2+\alpha}{2+\delta}(\ell -2), -\alpha)$.

    Using those estimates for the right hand side of \eqref{eq:formulaforstep2} and using the product rule $$[fg]_\alpha \leq \Vert f \Vert_\infty[g]_\alpha+[f]_\alpha \Vert g \Vert_\infty$$ on the term $\nabla \ln g_0 \otimes \nabla\ln g_0$, we get
    \begin{equation}
    \label{eq:alphacontrol}
    \left[\frac{\nabla^2 g_0}{g_0}\right]_{\alpha, \Sigma_R(0)} \leq C (R^\epsilon + R^{\beta+\theta}) \leq CR^{\beta + \theta}.
    \end{equation}

    Doing the same for the $L^\infty$ norm, the computations yield
    \begin{equation}
    \label{eq:infinitycontrol}
        \left\Vert\frac{\nabla^2 g_0}{g_0}\right\Vert_{L^\infty( \Sigma_R(0))} \leq C R^{2\beta}.
    \end{equation}

    \textit{Step 2.} Now that we now how $\frac{\nabla^2 g_0}{g_0}$ grows for large velocities, we work on the small cube $\Sigma_r(v)$. Once again by the product rule for Hölder norms,
    \begin{align*}
        [\nabla^2 g_0]_{\alpha,\Sigma_r(v)} &\leq \Vert g_0\Vert_{L^\infty( \Sigma_r(v))}\left[\frac{\nabla^2 g_0}{g_0}\right]_{\alpha, \Sigma_r(v)} + [g_0]_{\alpha, \Sigma_r(v)}\left\Vert\frac{\nabla^2 g_0}{g_0}\right\Vert_{L^\infty( \Sigma_r(v))}.
    \end{align*}
    Since $\Sigma_r(v) \subset \Sigma_{\vert v \vert+1}(0)$, we can use the two controls \eqref{eq:alphacontrol} and \eqref{eq:infinitycontrol} we derived above (with $R=\vert v \vert+1$):
    \begin{align}
    \label{eq:interpol}
        [\nabla^2 g_0]_{\alpha,\Sigma_r(v)}
        &\leq C\left( \Vert g_0\Vert_{L^\infty( \Sigma_r(v))} \jap{v}^{\beta+\theta}+ [g_0]_{\alpha, \Sigma_r(v)}\jap{v}^{2\beta}\right).
    \end{align}
    We are satisfied with the first term. We wish to control the $C^\alpha$ seminorm on the right by $\Vert g\Vert_{L^\infty( Q_r(v))}$. We once again use the logarithm of $g_0$. For any $w,z \in \Sigma_r(v)$,
    \begin{align*}
        \frac{\left \vert g_0(w) -g_0(z)\right \vert}{g_0(v)} = \left\vert \exp(\ln g_0 (w)-\ln g_0 (v)) - \exp(\ln g_0 (z)-\ln g_0 (v)) \right\vert.
    \end{align*}
    Recalling that $r=\frac{1}{2}\jap{v}^{-\beta}$, we have
    \begin{align*}
    \vert \ln g_0 (w)-\ln g_0 (v) \vert &\leq \Vert \nabla \ln g_0 \Vert_{L^\infty(\Sigma_r(v))} \vert w-v\vert\\
    &\leq C\jap{v}^\beta r\\
    &\leq C\jap{v}^\beta \jap{v}^{-\beta}\\
    &\leq C,
    \end{align*}
    where the control of the $L^\infty$ norm was obtained thanks to \eqref{eq:step2gradcontrol}. The cancellation above is why we choose this value of $r$. The same holds with $z$ in place of $w$. The arguments in the exponential are hence bounded, so by local Lipschitz behavior of the exponential, we have:
    \begin{align*}
         \left\vert \exp(\ln g_0 (w)-\ln g_0 (v)) - \exp(\ln g_0 (z)-\ln g_0 (v)) \right\vert &\leq C \vert \ln g_0 (w) - \ln g_0(z) \vert\\
         &\leq C [\ln g_0]_{\alpha, \Sigma_r(v)}\vert w-z\vert^\alpha.
    \end{align*}
    We can thus write
    \begin{align*}
    [g_0]_{\alpha, \Sigma_r(v)} &\leq C g_0(v) [\ln g_0]_{\alpha, \Sigma_r(v)}.
    \end{align*}
     Using that $[\ln g_0]_{\alpha, \Sigma_r(v)}\lesssim \jap{v}^{\sigma}$ with $\sigma=\max(2+\frac{\alpha}{2+\delta}(\ell-2),2-\alpha)$ (obtained by interpolation as above), we get
     \begin{align*}
    [g_0]_{\alpha, \Sigma_r(v)}
    &\leq C \Vert g_0\Vert_{L^\infty( \Sigma_r(v))} \jap{v}^{\sigma}.
    \end{align*}
     Plugging this back in \eqref{eq:interpol}, we reach
    \begin{align*}
        [\nabla^2 g_0]_{\alpha,\Sigma_r(v)} &\leq C \Vert g_0\Vert_{L^\infty( \Sigma_r(v))} \left(\jap{v}^{\beta+\theta}+ \jap{v}^{2\beta+\sigma}\right)\\
        &\leq C \Vert g_0\Vert_{L^\infty( \Sigma_r(v))}\jap{v}^{2\beta+\sigma},
    \end{align*}
    which is the control of the $C^{2+\alpha}$ norm we wanted. We rewrite it in terms of weighted norms for an easier use in the Schauder estimate: the weight on $\Sigma_r(v)$ is smaller than $r^{2+\alpha}$, so that
    \begin{equation}
    [g_0]^*_{2+\alpha,\Sigma_r(v)} \leq C r^{2+\alpha}\Vert g_0\Vert_{L^\infty( \Sigma_r(v))}\jap{v}^{2\beta+\sigma} = C\Vert g_0\Vert_{L^\infty( \Sigma_r(v))}\jap{v}^{-\beta\alpha+\sigma}.
    \end{equation}
    The proof is over.
\end{proof}

% Lemmas \ref{lem:Höldera} and \ref{lem:Hölderc} ensure that for some $\alpha>0$, the coefficients of the Landau equation in parabolic form \eqref{eq:landaupara} are $\alpha$-Hölder continuous in velocity. Because both coefficients are also uniformly bounded on $\mathbb{R}^3$, we can interpolate and obtain their $\beta$-Hölder continuity for any $0<\beta <\alpha$. Together with the ellipticity established in the previous section, this is what we needed to setup Schauder estimates. Indeed, we will not need any Hölder continuity in time of the coefficients $\bar{a}$ and $\bar{c}$. We wish to \textit{propagate} the Hölder continuity at initial time rather than \textit{generate} it to obtain uniform-in-time results that we can use in Theorem \ref{thm:rentropy}.

\appendix
\section{Maximum principle for the Landau equation}
\label{app:proofmaxprinc}
We state and prove the maximum principle we used to propagate the Maxwellian bounds (Proposition \ref{prop:propbounds}). Recall that the coefficients $\bar{a}$, $\bar{c}$ are defined in \eqref{eq:bara} and \eqref{eq:barc} with $g$ being the fixed solution of the Landau equation given by Proposition \ref{prop:solg}. Remark that the exponential growth condition below is largely satisfied in practice: we applied the principle to bounded functions.
\begin{prop}[Maximum principle]
\label{thm:maxprinciple}
    Let $\gamma \in [-3,-2]$. Let $u\in C([0,T]\times \mathbb{R}^3)$ be $C^1$ on $(0,T]$ and $C^2$ in space. Suppose that $u$ is a subsolution to
    $$\partial_t u \leq  \bar{a} : \nabla^2u  +\bar{c}u$$
    on $(0,T]\times \mathbb{R}^3$, which is non-positive at initial time:
    $$ \forall v \in \mathbb{R}^3\text{, }u(0,v) \leq 0.$$
    Further suppose that $u$ satisfies the growth condition
    $$ u(t,v) \leq  K e^{\eta \vert v \vert^2}$$
    for all $(t,v) \in [0,T]\times\mathbb{R}^3$, for some $K, \eta >0$. Then $u$ remains non-positive for all times, \textit{i.e.}
     $$ \forall (t,v) \in [0,T]\times\mathbb{R}^3\text{, }u(t,v) \leq 0.$$
\end{prop}

\begin{proof}
    We first reduce to the case of a non-positive zero order term. Recalling that $\bar{c}$ is bounded on $[0,T]\times\mathbb{R}^3$ (by Lemma \ref{lem:boundc}),
    the function $\tilde{u}(t,v)=e^{-\Vert \bar{c}\Vert_\infty t} u(t,v)$ satisfies
    $$\partial_t \tilde{u} \leq  \bar{a} : \nabla^2\tilde{u}  +\bar{c}\tilde{u} - \Vert \bar{c}\Vert_\infty \tilde{u} = \bar{a} : \nabla^2\tilde{u}  +\tilde{c}\tilde{u}$$
    with $\tilde{c} = \bar{c}-\Vert \bar{c}\Vert_\infty \leq 0$. Moreover,
    $\tilde{u}(t,v) \leq  K e^{-\Vert \bar{c}\Vert_\infty t}e^{\eta \vert v \vert^2} \leq K e^{\eta \vert v \vert^2}$ so the same growth condition holds.
    We now build a positive strict supersolution to the above equation. Let $\psi(t,v)=\exp(\omega t +2\eta \vert v \vert^2)>0$ for some $\omega>0$ to be determined.
    A direct computation yields:
    $$\nabla^2 \psi = (4\eta I + 16 \eta^2 v\otimes v)\psi.$$
    Using both directions in the ellipticity estimate (Lemma \ref{lem:elliptichigher}),
    $$\bar{a}:\nabla^2 \psi\ (v) \leq C(4 \eta \Lambda \jap{v}^{\gamma+2} + 16 \eta^2 \Lambda \jap{v}^{\gamma} \vert v \vert^2)\psi  \leq C\Lambda \eta(1+\eta)\jap{v}^{\gamma+2} \psi. $$
    Since $\gamma+2\leq 0$, this provides a uniform bound. Since $\tilde{c}\leq 0$,
    $$\bar{a} : \nabla^2 \psi  +\tilde{c} \psi \leq C\Lambda \eta(1+\eta)\psi.$$
    Using that $$\partial_t \psi(t,x)= \omega \psi,$$
    we can choose $\omega=C\Lambda \eta(1+\eta)+1$, such that $\psi$ is a strict supersolution:
    $$\partial_t \psi >  \bar{a} : \nabla^2 \psi  +\tilde{c} \psi.$$
    Let now $\varphi = u - \varepsilon \psi$ for some $\varepsilon>0$. By linearity, $\varphi$ is a strict subsolution, and using the growth bound, for all $(t,v) \in [0,T]\times\mathbb{R}^3$,
    $\varphi(t,v) \leq K e^{\eta \vert v \vert^2} - \varepsilon e^{2\eta \vert v \vert^2}.$
    Hence for $R$ large enough depending on $\varepsilon, \eta$ and $K$, $\varphi(t,v) < 0$ as soon as $\vert v \vert \geq R$. Now consider $(t_0, v_0)$ realizing the maximum of the continuous function $\varphi$ on the compact $[0,T]\times \overline{B(0,R)}$. Suppose that $\varphi(t_0,v_0)>0$, then necessarily $\vert v_0 \vert< R$, and $t_0>0$ because $\varphi$ is non-positive at initial time. Hence by the usual optimality conditions $\nabla^2 \varphi(t_0,v_0) \leq 0$, and $\partial_t \varphi(t_0,v_0)\geq 0$. But $\varphi$ being a strict subsolution,
    $$0 \leq \partial_t \varphi(t_0,v_0) < \bar{a}:\nabla^2 \varphi(t_0,v_0) +\tilde{c} \varphi(t_0,v_0) \leq \tilde{c} \varphi(t_0,v_0) \leq 0 $$
    since $\tilde{c} \leq 0$, which is absurd. Hence $\varphi(t_0,v_0)\leq 0$, so $u \leq \varepsilon \psi$ for all $(t,v) \in [0,T]\times\mathbb{R}^3$. Letting $\varepsilon \rightarrow 0$ concludes.
\end{proof}

\section{Proof of the Schauder estimates}
\label{app:schauder}
In this appendix we prove the Schauder estimate (Proposition \ref{thm:schauderlandau}). The proof is fairly standard but we write it in details to track the dependencies of the constants. It also differs a bit from the most classical setting because we work with regularity in velocity only. We begin with a Schauder estimate for constant-in-velocity pure-second-order differential operators, which is adapted from \cite{Brandt1969}.
\begin{prop}[Schauder estimate for frozen coefficients]
\label{prop:frozenschauder}
    Let $R>0$, and consider the cube $Q=Q_R(0)$.
    Let $L=\sum_{i,j}b_{ij}(t)\partial_{ij}$ be a uniformly elliptic operator with constant-in-$v$ coefficients, \textit{i.e.} there exists $\omega,\Omega>0$ such that, for any unit vector $\xi \in \mathbb{R}^3$, for all $t\in[0,T]$, the symmetric matrix $b(t)=(b_{ij}(t))$ satisfies
    $$\omega \leq \jap{\xi,b(t)\xi} \leq \Omega.$$
    Let $u:\bar{Q}\rightarrow \mathbb{R}$ be $C^1$ in time and $C^2$ in velocity (up to the boundary), and let $u_0:=u(0,\cdot)$.
    For any $t\in[0,T]$, any $\alpha \in (0,1)$, any directions $i,j,k\in\{1,2,3\}$, and any $h\in]0,R]$, it holds
    \begin{multline}
    \label{eq:schauderfrozen}
        \frac{\left\vert \partial_{ij}u(t,he_k)-\partial_{ij}u(t,-he_k)\right\vert}{(2h)^\alpha} \leq C_\alpha \big([u_0]_{2+\alpha,\Sigma_R}+\omega^{-1}[\partial_t u-Lu]_{\alpha,Q}\\
        +R^{-\alpha}(1+\omega^{-1}\Omega)\Vert \nabla^2 u \Vert_{L^\infty(Q)}\big)
    \end{multline}
    where $(e_l)_{l=1,2,3}$ is the canonic basis of $\mathbb{R}^3$ and $C_\alpha$ depends on $\alpha$ only.
\end{prop}
\begin{proof}
We adapt the proof in \cite{Brandt1969}, which is based on the maximum principle for $L$. In \cite{Brandt1969}, Hölder continuity is generated at the end time ($t=T$ in our notations) rather than propagated on $[0,T]$. We mostly keep the same notations apart from the new function $\psi_4$ which ensures this propagation.  We assume every quantity on the right hand side to be finite otherwise there is nothing to prove.
\\
\textit{Step 0.} We set up the framework of the proof. The key is to add a new independent variable $y=(y_1,y_2,y_3)\in\mathbb{R}^3$ and consider the $7$-dimensional domain of triplets $(t,v,y)$ in
$$\mathcal{Q}:=[0,T]\times(-R/4,R/4)^3\times(0,R/4)^3,$$
on which we define the uniformly elliptic operator $$\mathcal{L}:=L+\frac{\omega}{2} \Delta_y - \frac{\omega}{2} (\partial_{v_i v_i}+\partial_{v_j v_j}+\partial_{v_k v_k}).$$
We introduce a notation for the finite difference of a function $\varphi(v)$ in the direction $l$ with step $h$: $$\delta_{l}(h) \varphi(v) := \frac{1}{2}[\varphi(v+he_l)-\varphi(v-he_l)].$$
\\
\textit{Step 1.} Let define a function on $\mathcal{\bar{Q}}$ by: $$\phi(t,v,y) := \delta_k(y_1)\delta_j(y_2)\delta_i(y_3)u(t,v)$$
where the finite differences are iterated. One can note that the finite differences are well defined (the arguments of $u$ stay in $\bar{Q}$ thanks to the bounds on $v$ and $y$) and they commute. Using the mean-value theorem on any two of those, one easily gets:
\begin{equation}
    \label{eq:majphi}
    \vert \phi(t,v,y) \vert \leq C \Vert \nabla^2 u \Vert_{L^\infty(Q)} Y_{\min}
\end{equation}
with $Y_{\min}=\min(y_1 y_2, y_2 y_3 , y_3 y_1)$.

Using $\partial_{y_1 y_1} \varphi(v \pm y_1 e_k) = \partial_{v_k v_k} \varphi(v \pm y_1 e_k)$ and similar identities,
\begin{equation*}
    (\partial_t -\mathcal{L})\phi(t,v,y) = (\partial_t -L)\phi(t,v,y) = \delta_k(y_1)\delta_j(y_2)\delta_i(y_3)(\partial_t -L)u(t,v).
\end{equation*}
Here we have used that $L$ has constant coefficients to make it commute with the finite differences. Letting $H:=[\partial_t u-Lu]_{\alpha,Q}$, we have
\begin{equation}
    \label{eq:Lphi}
    \vert(\partial_t -\mathcal{L})\phi(t,v,y) \vert \leq H2^{\alpha-1} y_{\min}^{\alpha},
\end{equation}
with $y_{\min} = \min(y_1 ,y_2,y_3)$.

\textit{Step 2.} With those two bounds on $\phi$, we now wish to build an auxiliary function $\Psi$. We first define
$$\psi_1(t,v,y)=\psi_1(y)=y_1 y_2 y_3(y_1^\alpha+y_2^\alpha+y_3^\alpha)^{-\frac{1-\alpha}{\alpha}} \geq 0.$$
A lengthy but straightforward computation yields
$$\mathcal{L}\psi_1 = \frac{\omega}{2}\Delta \psi \leq -C_\alpha \omega (y_1^\alpha+y_2^\alpha+y_3^\alpha)^{-\frac{1+\alpha}{\alpha}} \sum_{l\neq m\neq n} y_l^{\alpha -1} y_m y_n^{\alpha +1} .$$
Selecting a term of the sum with $l,m,n$ such that $y_{\min} = y_l \leq y_m \leq y_n$,
\begin{align*}
    \mathcal{L}\psi_1 &\leq -C_\alpha \omega (3y_n^\alpha)^{-\frac{1+\alpha}{\alpha}} y_l^{\alpha -1} y_m y_n^{\alpha +1}\\
    &\leq -C_\alpha \omega y_l^{\alpha-1}y_l\\
    &\leq -C_\alpha \omega y_{\min}^\alpha.
\end{align*}
Hence as soon as $C_\alpha>0$ is large enough, 
$$(\partial_t -\mathcal{L})(C_\alpha H\omega^{-1} \psi_1) \geq H2^{\alpha-1} y_{\min}^\alpha \geq \vert(\partial_t -\mathcal{L})\phi(t,v,y) \vert, $$
recalling \eqref{eq:Lphi}.

We now define
$$\psi_2(t,v,y)=\psi_2(y)=\frac{4}{R}y_1 y_2 y_3 \geq 0,$$
which satisfies
$$(\partial_t - \mathcal{L})\psi_2=0\geq 0,$$
and if $y_l=\frac{R}{4}$ for some $l$, $\psi_2\geq Y_{\min}$, so recalling \eqref{eq:majphi},
$$C \Vert\nabla^2u\Vert_{L^\infty(Q)}\psi_2 \geq\vert \phi\vert $$
on the $\{y_l=\frac{R}{4} \}$ part of the boundary of $\mathcal{Q}$, again for $C$ large enough. On the $\{y_l=0 \}$ part of the boundary, $\vert \phi \vert = 0 \leq \psi_2$.

For the boundary in $v$, we define
$$\psi_3(t,v,y)=\psi_3(v,y) = \frac{16}{R^2}y_2 y_3 \left(\sum_{i=1}^3 v_i^2 +\frac{9\Omega}{\omega}y_1\left(\frac{R}{4}-y_1\right)\right)\geq 0,$$
and we have
$$(\partial_t - \mathcal{L})\psi_3 = \frac{16}{R^2}\left[-y_2 y_3 \sum_{i=1}^3 2b_{ii} + 9y_2 y_3 \Omega - 3 \omega y_2 y_3\right] \geq \frac{16}{R^2}y_2 y_3 (-6\Omega +9\Omega -3 \Omega)\geq 0$$
and if $\vert v_l\vert=\frac{R}{4}$ for some $l$, we have $\psi_3 \geq y_2 y_3 \geq Y_{\min}$, so, again for $C$ large enough,
$$C \Vert\nabla^2u\Vert_{L^\infty(Q)}\psi_3 \geq\vert \phi\vert $$
on the $\{\vert x_l\vert=\frac{R}{4} \}$ part of the boundary.

It remains to treat the boundary in time. We let 
$$\psi_4(t,v,y)=\psi_4(y)=y_1^\alpha y_2 y_3\geq 0,$$
which satisfies
$$(\partial_t - \mathcal{L})\psi_4 = -\frac{\omega}{2} \alpha (\alpha-1)y_1^{\alpha -2}y_2 y_3 \geq 0,$$
and at $t=0$, iterating two times the mean value theorem, and using Hölder continuity of $\nabla^2u_0$,
$$\vert \phi(0,v,y)\vert \leq C [\nabla^2 u_0]_{\alpha,(-R,R)^3} y_1^\alpha y_2 y_3 = C [u_0]_{2+\alpha,(-R,R)^3} \psi_4(y).$$
Hence for $C_\alpha$ large enough, gathering every estimate, the function
$$\Psi := C_\alpha \left[ H\omega^{-1} \psi_1 + \Vert\nabla^2u\Vert_{L^\infty(Q)}(\psi_2 + \psi_3) +[u_0]_{2+\alpha,(-R,R)^3}\psi_4\right]$$
satisfies
$$(\partial_t-\mathcal{L})\Psi \geq \vert(\partial_t -\mathcal{L})\phi \vert $$
and
$$\Psi \geq \vert \phi \vert \, \, \text{ on } \bar{\mathcal{Q}}\setminus \left((0,T]\times(-R/4,R/4)^3\times(0,R/4)^3 \right).$$
\\
\textit{Step 4.} Applying the maximum principle  to $\pm \phi-\Psi$, we get $\pm \phi-\Psi\leq 0$, \textit{i.e.} $\vert \phi \vert \leq \Psi$, on $\mathcal{Q}$. Evaluating at $v=0$, it yields
\begin{align*}
    \vert \phi(t,0,y)\vert  \leq &C_\alpha [H \omega^{-1} y_1 y_2 y_3(y_1^\alpha+y_2^\alpha+y_3^\alpha)^{-\frac{1-\alpha}{\alpha}}  \\
    &+ \Vert\nabla^2u\Vert_{L^\infty(Q)}\left(\frac{4}{R}y_1 y_2 y_3 +  \frac{16}{R^2}y_1y_2 y_3 \frac{9\Omega}{\omega}\left(\frac{R}{4}-y_1\right)\right)\\
    &+[u_0]_{2+\alpha,(-R,R)^3}y_1^\alpha y_2 y_3].
\end{align*}
Dividing by $y_1^\alpha y_2 y_3$, plugging in the bounds
$$y_1^{1-\alpha}(y_1^\alpha+y_2^\alpha+y_3^\alpha)^{-\frac{1-\alpha}{\alpha}} \leq 1$$
$$\frac{y_1^{1-\alpha}}{R}\leq R^{-\alpha}$$
in the terms that appear and dropping the only negative term,
\begin{align*}
    \frac{\vert \phi(t,0,y)\vert}{y_1^\alpha y_2 y_3}  \leq &C_\alpha \left[ H\omega^{-1}
    + \Vert\nabla^2u\Vert_{L^\infty(Q)}\left(R^{-\alpha} \left(1+\frac{\Omega}{\omega}\right)\right)+[u_0]_{2+\alpha,(-R,R)^3}\right],
\end{align*}
up to a change of $C_\alpha$. The right hand side is the same as in \eqref{eq:schauderfrozen}. Letting $y_2,y_3 \rightarrow 0$, and taking $y_1 = h$, we get
\begin{align*}
    &\frac{\vert \partial_{ij}u(t,he_k) - \partial_{ij}u(t,-he_k)\vert}{(2h)^\alpha} \\
    &\;\;\;\;\leq C_\alpha \left[ H\omega^{-1}
    + \Vert\nabla^2u\Vert_{L^\infty(Q)}\left(R^{-\alpha} \left(1+\frac{\Omega}{\omega}\right)\right)+[u_0]_{2+\alpha,(-R,R)^3}\right],
\end{align*} 
which is \eqref{eq:schauderfrozen} for $h\leq\frac{R}{4}$. But if $h\geq \frac{R}{4}$, \eqref{eq:schauderfrozen} is a direct consequence of
$$
    \frac{\left\vert \partial_{ij}u(t,he_k)-\partial_{ij}u(t,-he_k)\right\vert}{(2h)^\alpha} \leq C_\alpha R^{-\alpha}\Vert \nabla^2 u \Vert_{L^\infty(Q)}.
$$
\\
\textit{Bonus Step.} The proof of the maximum principle we used is as follows. Let $$\zeta:=\phi-\Psi - \varepsilon(1+t)$$ for some small $\varepsilon>0$. We have
$$(\partial_t-\mathcal{L})\zeta\leq -\varepsilon<0  \, \, \text{ on }\bar{\mathcal{Q}}$$
and
$$\zeta \leq 0  \, \, \text{ on } \bar{\mathcal{Q}}\setminus \left((0,T]\times(-R/4,R/4)^3\times(0,R/4)^3 \right).$$
Suppose the maximum of the continuous function $\zeta$ on $\bar{\mathcal{Q}}$ is achieved at a point $$(t_0,v_0,y_0)\in(0,T]\times(-R/4,R/4)^3\times(0,R/4)^3.$$
Then $\nabla^2_{v,y} \zeta(t_0,v_0,y_0) \leq 0$, $\partial_t\zeta(t_0,v_0,y_0) \geq 0$, hence
$$(\partial_t-\mathcal{L})\zeta(t_0,v_0,y_0) \geq 0$$ which is absurd. So $\zeta\leq 0$ on $\bar{\mathcal{Q}}$ and letting $\varepsilon\rightarrow 0$ concludes. The same goes for the function $-\phi-\Psi$.
\end{proof}

Using Proposition \ref{prop:frozenschauder}, we can obtain the Schauder estimate (Proposition \ref{thm:schauderlandau}) by perturbation of the coefficients.

\begin{proof}\textit{(Proposition \ref{thm:schauderlandau})}

        We assume the right hand side of \eqref{eq:schauderlandau} to be finite.
        It is easily verified that Hölder continuity can be checked coordinate by coordinate, \textit{i.e.} the seminorm
        $$[[\varphi]]_{\alpha,Q_{R}(v_0)} := \sup_{k=1,2,3}\sup_{(t,v)\in Q_{R}}\sup_{\substack{h\in \mathbb{R^*} \\(t,v+he_k)\in Q_{R}}} \frac{\vert \varphi(t,v)-\varphi(t,v+he_k)\vert}{\vert h\vert^\alpha},$$
        is equivalent to $[\varphi]_{\alpha,Q_{R}(v_0)}$. This transposes to the weighted version of these norms. In particular, for some absolute constant $C_{eq}>0$, by definition of the supremum, there exists $i,j,k\in\{1,2,3\}$, $(t,v_2)\in Q_{R}(v_0)$, $\tilde{h}\in \mathbb{R^*}$, so that
        $$C_{eq}[u]^*_{2+\alpha,Q_{R}(v_0)} \leq d_{v_2,v_2+\tilde{h}e_k}^{2+\alpha} \frac{\vert \partial_{ij}u(t,v_2)- \partial_{ij}u(t,v_2+\tilde{h}e_k)\vert}{\vert \tilde{h}\vert^\alpha}.$$
        Considering the middle point $v_1=v_2+\frac{\tilde{h}}{2}e_k$, letting $d_1=d_{v_2,v_2+\tilde{h}e_k}$ and $h=\vert\frac{\tilde{h}}{2}\vert$, we can rewrite this as:
        $$C_{eq}[u]^*_{2+\alpha,Q_{R}(v_0)} \leq d_1^{2+\alpha} \frac{\vert \partial_{ij}u(t,v_1+he_k)- \partial_{ij}u(t,v_1-he_k)\vert}{( 2h)^\alpha}.$$
        Let $0<\eta<\frac{d_1}{2}$ to be determined later.
        If $h\geq \eta$,
        \begin{align}
        \label{eq:bigh}
            C_{eq}[u]^*_{2+\alpha,Q_{R}(v_0)} &\leq d_1^{2+\alpha} (2\eta)^{-\alpha}\left(\vert \partial_{ij}u(t,v_1+he_k)\vert + \vert \partial_{ij}u(t,v_1-he_k)\vert\right)\nonumber\\
            &\leq 2\ d_1^\alpha \ (2\eta)^{-\alpha} \Vert d^2 \nabla^2 u\Vert_{L^\infty(Q_R(v_0))}
         \end{align}
        since by definition $d_1 = \min(d_{v_1+h e_k},d_{v_1-h e_k})$. We keep this bound in mind for later.
        
        Otherwise, $h<\eta$, so $v_1\pm h e_k \in Q_{\eta}(v_1)\subset  \overline{Q_{\eta}(v_1)}\subset Q_R(v_0)$. This last inclusion is true thanks to the condition $\eta<\frac{d_1}{2} $, remembering that $d_1$ measures the distance of $v_1\pm h e_k$ to the edge of $\Sigma_R(v_0)$.
        It ensures that $u$ is $C^1$ in $t$, $C^2$ in $v$, up to the boundary of $\overline{Q_{\eta}(v_1)}$, so that we can apply Proposition \ref{prop:frozenschauder}. With this in mind, we define the frozen-in-$v$ operator
        $$L = \bar{a}(t,v_1):\nabla^2$$
        which, by Lemmas \ref{lem:ellipticlower} and \ref{lem:elliptichigher}, is uniformly elliptic with constants depending on $v_1$:
        $$\omega := \lambda\jap{v_1}^\gamma \leq \jap{\bar{a}(t,v_1)\xi,\xi} \leq \Lambda \jap{v_1}^{\gamma+2} =: \Omega.$$
        We apply Proposition \ref{prop:frozenschauder} on $Q_\eta(v_1)$ (in all rigor, we apply it to the translated function $u(t,v_1+v)$ in $Q_\eta(0)$), yielding:
        \begin{align*}
        d_1^{2+\alpha}\frac{\left\vert \partial_{ij}u(t,v_1+he_k)-\partial_{ij}u(t,v_1-he_k)\right\vert}{(2h)^\alpha} \leq d_1^{2+\alpha}C_\alpha &([u_0]_{2+\alpha,\Sigma_\eta(v_1)}+\omega^{-1}[\partial_t u-Lu]_{\alpha,Q_\eta(v_1)}\\
        &+\eta^{-\alpha}(1+\omega^{-1}\Omega)\Vert \nabla^2 u \Vert_{L^\infty(Q_{\eta}(v_1))}).
        \end{align*}
        For all $v\in Q_\eta(v_1)$, $d_v\geq \frac{d_1}{2}$ (once again thanks to $\eta<\frac{d_1}{2}$). So
        $$d_1^{2+\alpha}[u_0]_{2+\alpha,\Sigma_\eta(v_1)}\leq 2^{2+\alpha} [u_0]^*_{2+\alpha,\Sigma_R(v_0)},$$
        and similarly
        $$d_1^{2+\alpha}\eta^{-\alpha}\Vert \nabla^2 u \Vert_{L^\infty(Q_{\eta}(v_1))}\leq 4d_1^\alpha  \eta^{-\alpha} \Vert d^2 \nabla^2 u \Vert_{L^\infty(Q_{R}(v_0))}.$$
        We now have to deal with the term
        \begin{align*}
            (\partial_t - L)u = \left(\bar{a}(t,v)-\bar{a}(t,v_1)\right):\nabla^2 u + \bar{c}u.
        \end{align*}
        Using $[\varphi \psi]_\alpha \leq [\varphi]_\alpha \Vert \psi\Vert_{L^\infty} +\Vert \varphi\Vert_{L^\infty}[\psi]_\alpha$,
        we get
        \begin{align*}
            d_1^{2+\alpha}[\bar{c}u]_{\alpha,Q_\eta(v_1)} &\leq d_1^{2+\alpha}[\bar{c}]_{\alpha,Q_\eta(v_1)}\Vert u\Vert_{L^\infty(Q_\eta(v_1))} +d_1^{2+\alpha}\Vert \bar{c}\Vert_{L^\infty(Q_\eta(v_1))}[u]_{\alpha,Q_\eta(v_1)}\\
            &\leq 2^{2+\alpha}\left([\bar{c}]_\alpha^\bigstar\Vert u\Vert_{L^\infty(Q_R(v_0))} +\Vert d^2 \bar{c}\Vert_{L^\infty(Q_R(v_0))}[u]^*_{\alpha,Q_R(v_0)}\right),
        \end{align*}
        using once again that on $Q_\eta(v_1)$, $d_v\geq \frac{d_1}{2}$, and recalling that $[\bar{c}]_\alpha^\bigstar$ is specifically defined with a $2+\alpha$ weight. 
        Then,
        \begin{align*}
        [\left(\bar{a}-\bar{a}(\cdot,v_1)\right):&\nabla^2 u ]_{\alpha,Q_\eta(v_1)} \\
        &\leq[\bar{a}]_{\alpha,Q_\eta(v_1)}\Vert \nabla^2 u\Vert_{L^\infty(Q_\eta(v_1))} +\Vert \bar{a}-\bar{a}(\cdot,v_1)\Vert_{L^\infty(Q_\eta(v_1))}[u]_{2+\alpha,Q_\eta(v_1)},
        \end{align*}
        where we have used that the constant-in-$v$ term $\bar{a}(\cdot,v_1)$ disappears from the seminorm, so that
        \begin{align*}
            d_1^{2+\alpha}[&\left(\bar{a}-\bar{a}(\cdot,v_1)\right):\nabla^2 u ]_{\alpha,Q_\eta(v_1)}\\
            &\leq 2^{2+\alpha} \left( [\bar{a}]^*_{\alpha,Q_R(v_0)}\Vert d^2\nabla^2 u\Vert_{L^\infty(Q_R(v_0))} +\Vert \bar{a}-\bar{a}(\cdot,v_1)\Vert_{L^\infty(Q_\eta(v_1))}[u]^*_{2+\alpha,Q_R(v_0)}\right),
        \end{align*}
        by adequately sharing the weight.
        Now, the Hölder continuity of $\bar{a}$ allows to write:
        $$\Vert \bar{a}-\bar{a}(\cdot,v_1)\Vert_{L^\infty(Q_\eta(v_1))}\leq C [\bar{a}]_{\alpha,Q_\eta(v_1)}\eta^\alpha \leq C 2^\alpha [\bar{a}]^*_{\alpha,Q_R(v_0)}\left(\frac{\eta}{d_1}\right)^\alpha.$$
        We have finally expressed everything in terms of weighted seminorms. Gathering everything, up to a change in $C_\alpha$,
        \begin{align*}
            C_\alpha^{-1}C_{eq}&[u]^*_{2+\alpha,Q_{R(v_0)}}\\
            \leq  &\ [u_0]^*_{2+\alpha,\Sigma_R(v_0)}\\
            &+(1+\omega^{-1}\Omega)d_1^\alpha  \eta^{-\alpha} \Vert d^2 \nabla^2 u \Vert_{L^\infty(Q_{R}(v_0))}\\
            &+\omega^{-1}\left([\bar{c}]_\alpha^\bigstar\Vert u\Vert_{L^\infty(Q_R(v_0))} +\Vert d^2 \bar{c}\Vert_{L^\infty(Q_R(v_0))}[u]^*_{\alpha,Q_R(v_0)}\right)\\
            &+\omega^{-1}\left( [\bar{a}]^*_{\alpha,Q_R(v_0)}\Vert d^2\nabla^2 u\Vert_{L^\infty(Q_R(v_0))} +[\bar{a}]^*_{\alpha,Q_R(v_0)}\left(\frac{\eta}{d_1}\right)^\alpha[u]^*_{2+\alpha,Q_R(v_0)}\right)
        \end{align*}
        Remembering the loose end \eqref{eq:bigh} that was for $h\geq \eta$, we see that \eqref{eq:bigh} implies the above bound, so we are now back to the general case.
        We choose $\eta$ small enough, such that
        $$\omega^{-1}[\bar{a}]^*_{\alpha,Q_R(v_0)}\left(\frac{\eta}{d_1}\right)^\alpha = \frac{1}{2}C_\alpha^{-1}C_{eq},$$
        which allows to pass to the left-hand-side the $2+\alpha$ weighted norm of $u$ (the very last term in the above inequality).
        We also have $d_1^\alpha  \eta^{-\alpha} \approx \omega^{-1}[\bar{a}]^*_{\alpha,Q_R(v_0)}$, so, changing $C_\alpha$ again, it holds that
        \begin{align*}
            C_\alpha^{-1}[u]^*_{2+\alpha,Q_{R(v_0)}} \leq  &\ [u_0]^*_{2+\alpha,\Sigma_R(v_0)}\\
            &+\omega^{-1}(1+\omega^{-1}\Omega)[\bar{a}]^*_{\alpha,Q_R(v_0)} \Vert d^2 \nabla^2 u \Vert_{L^\infty(Q_{R}(v_0))}\\
            &+\omega^{-1}\left([\bar{c}]_\alpha^\bigstar\Vert u\Vert_{L^\infty(Q_R(v_0))} +\Vert d^2 \bar{c}\Vert_{L^\infty(Q_R(v_0))}[u]^*_{\alpha,Q_R(v_0)}\right)\\.
        \end{align*}
        We use interpolation to get rid of derivatives of $u$: by chaining the interpolations from Lemma \ref{lem:interpolation2}, for $K,K'\geq 1$, 
         $$K\Vert d^2 \nabla^2 u \Vert_{L^\infty(Q_{R}(v_0))}\leq \frac{1}{4}[u]^*_{2+\alpha,Q_{R(v_0)}} + C K^{1+\frac{2}{\alpha}}\Vert  u \Vert_{L^\infty(Q_{R}(v_0))}$$
         and
         $$K'[u]^*_{\alpha,Q_R(v_0)}\leq \frac{1}{4}[u]^*_{2+\alpha,Q_{R(v_0)}} + C K'^{1+\frac{\alpha}{2}}\Vert  u \Vert_{L^\infty(Q_{R}(v_0))}.$$
        Using those with $K= 1 + C_\alpha \omega^{-1}(1+\omega^{-1}\Omega)[\bar{a}]^*_{\alpha,Q_R(v_0)}$ and $K' = 1 + C_\alpha \omega^{-1} \Vert d^2 \bar{c}\Vert_{L^\infty(Q_R(v_0))}$, we finally obtain
        \begin{align*}
            C_\alpha^{-1}[u]^*_{2+\alpha,Q_{R(v_0)}} \leq  &\ [u_0]^*_{2+\alpha,\Sigma_R(v_0)}\\
            &+\left[1+\omega^{-1}(1+\omega^{-1}\Omega)[\bar{a}]^*_{\alpha,Q_R(v_0)}\right]^{1+\frac{2}{\alpha}} \Vert u \Vert_{L^\infty(Q_{R}(v_0))}\\
            &+\omega^{-1}[\bar{c}]_\alpha^\bigstar\Vert u\Vert_{L^\infty(Q_R(v_0))}\\
            &+\left[1+\omega^{-1}\Vert d^2 \bar{c}\Vert_{L^\infty(Q_R(v_0))}\right]^{1+\frac{\alpha}{2}} \Vert u\Vert_{L^\infty(Q_R(v_0))}.
        \end{align*}
        Replacing $\omega$ and $\Omega$ by their expression and using $\jap{v_1}\approx\jap{v_0}$
        (which holds because $R\leq \jap{v_0}/2$), we obtain \eqref{eq:schauderlandau}.
\end{proof}

\section{Interpolations for Hölder norms}
\label{app:interpol}
We recall the following interpolation result, which holds for unweighted norms:
\begin{lemma}
\label{lem:interpolation1}
    For any $\alpha\in(0,1)$, any function $u$ defined on a cube $Q=\Sigma_R(v_0)$ or $Q=Q_R(v_0)$ and any $0<\varepsilon <  R$, any $\beta\in(\alpha,1)$,
    \begin{align*}
        \Vert  \nabla u \Vert_{L^\infty(Q)}&\leq \varepsilon^\alpha [u]_{1+\alpha,Q} + C \varepsilon^{-1} \Vert  u \Vert_{L^\infty(Q)} \\
        \Vert  \nabla u \Vert_{L^\infty(Q)}&\leq \varepsilon\Vert\nabla^2 u\Vert_{L^\infty(Q)} + C \varepsilon^{-1} \Vert  u \Vert_{L^\infty(Q)} \\
        [ u ]_{\alpha,Q}&\leq \varepsilon^{\beta-\alpha}  [ u ]_{\beta,Q} + C \varepsilon^{-\alpha} \Vert  u \Vert_{L^\infty(Q)}, \\
        [ u ]_{\alpha,Q}&\leq \varepsilon^{1-\alpha} \Vert  \nabla u \Vert_{L^\infty(Q)} + C \varepsilon^{-\alpha} \Vert  u \Vert_{L^\infty(Q)}, \\
    \end{align*}
    for some $C>0$ independent of $u$ and $\varepsilon$.
\end{lemma}
% \begin{proof}
%     We can suppose that there is no time dependence, \textit{i.e.} $Q=\Sigma_R(v_0)$.
%     Let $v\in Q$. Because $\varepsilon<R$, for any $i\in\{1,2,3\}$, either $v+\varepsilon e_i$ or $v-\varepsilon e_i$ is in $Q$. Let denote $w=v+he_i$ such a point. By the mean value theorem,
%     there exists $z\in[v,w]$ such that
%     $$\partial_i u (z) = \frac{u(w)-u(v)}{h}.$$
%     Hence
%     \begin{align*}
%         \vert \partial_i u (v) \vert &\leq \vert \partial_i u (v) - \partial_i u (z) \vert +  \frac{\vert u(w)-u(v) \vert}{\varepsilon}\\
%         &\leq \vert v-z\vert^\alpha[\partial_i u]_{\alpha,Q} + 2\varepsilon^{-1}\Vert u \Vert_{L^\infty(Q)}\\
%         &\leq \varepsilon^\alpha [\partial_i u]_{\alpha,Q} + 2\varepsilon^{-1}\Vert u \Vert_{L^\infty(Q)}.
%     \end{align*}
%     This yields the first interpolation, and the second by using
%     $$\vert \partial_i u (v) - \partial_i u (z) \vert \leq \vert v-z\vert \Vert  \nabla^2 u \Vert_{L^\infty(Q)} $$ in the first line instead of Hölder continuity.

%     For the third , for $v,w\in Q$, if $\vert v-w\vert\leq \varepsilon$,
%     $$\vert u(v) - u(w) \vert \leq \vert v-w \vert^\alpha \varepsilon^{\beta-\alpha}[ u ]_{\beta,Q},$$
%     else $\vert v-w\vert> \varepsilon$ and
%     $$\vert u(v) - u(w) \vert \leq \vert v-w\vert^\alpha \varepsilon^{-\alpha} 2 \Vert u \Vert_{L^\infty(Q)}.$$
%     The fourth is proven similarly.
% \end{proof}
Those results can be chained together to control any norm between $C^0$ and $C^{2+\alpha}$ by those endpoints, with the $\varepsilon$ powers scaling accordingly. We also have a similar interpolation result for weighted norms, with no assumption on the size of the cube.
\begin{lemma}
\label{lem:interpolation2}
    For any $\alpha\in(0,1)$, any $C^{2+\alpha}$ function $u$ defined on a cube $Q=\Sigma_R(v_0)$ or $Q=Q_R(v_0)$, and any $0<\varepsilon <  \frac{1}{2}$,
    \begin{align*}
        \Vert  d^2\nabla^2 u \Vert_{L^\infty(Q)}&\leq C\left(\varepsilon^\alpha [ u]^*_{2+\alpha,Q} +  \varepsilon^{-1} \Vert d \nabla u \Vert_{L^\infty(Q)}\right) \\
        \Vert  d\nabla u \Vert_{L^\infty(Q)}&\leq C\left(\varepsilon^\alpha [ u]^*_{1+\alpha,Q} +  \varepsilon^{-1} \Vert  u \Vert_{L^\infty(Q)}\right) \\
        \Vert  d\nabla u \Vert_{L^\infty(Q)}&\leq C\left(\varepsilon\Vert d^2 \nabla^2 u\Vert_{L^\infty(Q)} +  \varepsilon^{-1} \Vert  u \Vert_{L^\infty(Q)}\right) \\
        [ u ]^*_{\alpha,Q}&\leq C\left( \varepsilon^{1-\alpha} \Vert d \nabla u \Vert_{L^\infty(Q)} +  \varepsilon^{-\alpha} \Vert  u \Vert_{L^\infty(Q)}\right), \\
    \end{align*}
    for some $C>0$ depending only on $\alpha$.
    
\end{lemma}

For proofs of these results and variants, we refer to textbooks in elliptic PDEs such as \cite[Chapter 3.2]{Krylov1996}, \cite[Chapter 6.8]{GilbargTrudinger2001}.

\sloppy
\printbibliography[title = References]

@article{AlonsoBaglandDesvillettes2024,
  title = {A {{Priori Estimates}} for {{Solutions}} to {{Landau Equation Under Prodi}}--{{Serrin Like Criteria}}},
  author = {Alonso, R. and Bagland, V. and Desvillettes, L. and Lods, B.},
  year = {2024},
  month = may,
  journal = {Archive for Rational Mechanics and Analysis},
  volume = {248},
  number = {3},
  pages = {42},
  issn = {1432-0673},
  doi = {10.1007/s00205-024-01992-y},
  urldate = {2025-02-28},
  langid = {english}
}

@article{Brandt1969,
  title = {Interior Schauder Estimates for Parabolic Differential- (or Difference-) Equations via the Maximum Principle},
  author = {Brandt, A.},
  year = {1969},
  month = sep,
  journal = {Israel Journal of Mathematics},
  volume = {7},
  number = {3},
  pages = {254--262},
  issn = {1565-8511},
  doi = {10.1007/BF02787619},
  urldate = {2025-01-21},
  langid = {english}
}

@article{CameronSilvestreSnelson2018a,
  title = {Global a Priori Estimates for the Inhomogeneous {{Landau}} Equation with Moderately Soft Potentials},
  author = {Cameron, Stephen and Silvestre, Luis and Snelson, Stanley},
  year = {2018},
  month = may,
  journal = {Annales de l'Institut Henri Poincar{\'e} C, Analyse non lin{\'e}aire},
  volume = {35},
  number = {3},
  pages = {625--642},
  issn = {0294-1449},
  doi = {10.1016/j.anihpc.2017.07.001},
  urldate = {2025-02-24}
}

@misc{CarrapatosoDesvillettesHe2015,
  title = {Estimates for the Large Time Behavior of the {{Landau}} Equation in the {{Coulomb}} Case},
  author = {Carrapatoso, Kleber and Desvillettes, Laurent and He, Lingbing},
  year = {2015},
  month = oct,
  number = {arXiv:1510.08704},
  eprint = {1510.08704},
  primaryclass = {math},
  publisher = {arXiv},
  doi = {10.48550/arXiv.1510.08704},
  urldate = {2025-01-07},
  archiveprefix = {arXiv}
}

@misc{CarrilloFengGuo2024,
  title = {Relative {{Entropy Method}} for {{Particle Approximation}} of the {{Landau Equation}} for {{Maxwellian Molecules}}},
  author = {Carrillo, Jos{\'e} Antonio and Feng, Xuanrui and Guo, Shuchen and Jabin, Pierre-Emmanuel and Wang, Zhenfu},
  year = {2024},
  month = sep,
  number = {arXiv:2408.15035},
  eprint = {2408.15035},
  primaryclass = {math},
  publisher = {arXiv},
  doi = {10.48550/arXiv.2408.15035},
  urldate = {2024-10-04},
  archiveprefix = {arXiv}
}

@misc{ChernGualdani2022,
  title = {Uniqueness of Higher Integrable Solution to the {{Landau}} Equation with {{Coulomb}} Interactions},
  author = {Chern, Jann-Long and Gualdani, Maria},
  year = {2022},
  month = jul,
  number = {arXiv:1910.06216},
  eprint = {1910.06216},
  publisher = {arXiv},
  urldate = {2024-11-07},
  archiveprefix = {arXiv}
}

@article{Dafermos1979,
  title = {The Second Law of Thermodynamics and Stability},
  author = {Dafermos, C. M.},
  year = {1979},
  month = jun,
  journal = {Archive for Rational Mechanics and Analysis},
  volume = {70},
  number = {2},
  pages = {167--179},
  issn = {1432-0673},
  doi = {10.1007/BF00250353},
  urldate = {2025-02-24},
  langid = {english}
}

@misc{Desvillettes2014a,
  title = {Entropy Dissipation Estimates for the {{Landau}} Equation in the {{Coulomb}} Case and Applications},
  author = {Desvillettes, Laurent},
  year = {2014},
  month = oct,
  number = {arXiv:1408.6025},
  eprint = {1408.6025},
  primaryclass = {math},
  publisher = {arXiv},
  doi = {10.48550/arXiv.1408.6025},
  urldate = {2025-02-11},
  archiveprefix = {arXiv}
}

@misc{DesvillettesGoldingGualdani2024,
  title = {Production of the {{Fisher}} Information for the {{Landau-Coulomb}} Equation with {{L1}} Initial Data},
  author = {Desvillettes, Laurent and Golding, William and Gualdani, Maria Pia and Loher, Amelie},
  year = {2024},
  month = oct,
  number = {arXiv:2410.10765},
  eprint = {2410.10765},
  primaryclass = {math},
  publisher = {arXiv},
  doi = {10.48550/arXiv.2410.10765},
  urldate = {2024-12-09},
  archiveprefix = {arXiv}
}

@article{Diperna1979,
  title = {Uniqueness of {{Solutions}} to {{Hyperbolic Conservation Laws}}},
  author = {Diperna, Ronald J.},
  year = {1979},
  journal = {Indiana University Mathematics Journal},
  volume = {28},
  number = {1},
  eprint = {24892290},
  eprinttype = {jstor},
  pages = {137--188},
  publisher = {Indiana University Mathematics Department},
  issn = {0022-2518},
  urldate = {2025-02-24}
}

@misc{Fournier2008,
  title = {Particle Approximation of Some {{Landau}} Equations},
  author = {Fournier, Nicolas},
  year = {2008},
  month = nov,
  number = {arXiv:0811.2688},
  eprint = {0811.2688},
  publisher = {arXiv},
  doi = {10.48550/arXiv.0811.2688},
  urldate = {2024-10-25},
  archiveprefix = {arXiv}
}

@article{Fournier2010,
  title = {Uniqueness of Bounded Solutions for the Homogeneous {{Landau}} Equation with a {{Coulomb}} Potential},
  author = {Fournier, Nicolas},
  year = {2010},
  month = nov,
  journal = {Communications in Mathematical Physics},
  volume = {299},
  number = {3},
  eprint = {0909.0647},
  primaryclass = {math},
  pages = {765--782},
  issn = {0010-3616, 1432-0916},
  doi = {10.1007/s00220-010-1113-9},
  urldate = {2025-01-29},
  archiveprefix = {arXiv}
}

@misc{FournierHauray2015,
  title = {Propagation of Chaos for the {{Landau}} Equation with Moderately Soft Potentials},
  author = {Fournier, Nicolas and Hauray, Maxime},
  year = {2015},
  month = jan,
  number = {arXiv:1501.01802},
  eprint = {1501.01802},
  primaryclass = {math},
  publisher = {arXiv},
  urldate = {2024-09-05},
  archiveprefix = {arXiv},
  langid = {english}
}

@book{GilbargTrudinger2001,
  title = {Elliptic {{Partial Differential Equations}} of {{Second Order}}},
  author = {Gilbarg, David and Trudinger, Neil S.},
  year = {2001},
  month = jan,
  publisher = {Springer Science \& Business Media},
  googlebooks = {eoiGTf4cmhwC},
  isbn = {978-3-540-41160-4},
  langid = {english}
}

@misc{GolseImbertJi2024,
  title = {Local Regularity for the Space-Homogeneous {{Landau}} Equation with Very Soft Potentials},
  author = {Golse, Fran{\c c}ois and Imbert, Cyril and Ji, Sehyun and Vasseur, Alexis F.},
  year = {2024},
  month = jan,
  number = {arXiv:2206.05155},
  eprint = {2206.05155},
  primaryclass = {math},
  publisher = {arXiv},
  doi = {10.48550/arXiv.2206.05155},
  urldate = {2025-02-12},
  archiveprefix = {arXiv}
}

@misc{GuerinFournier2008,
  title = {Well-Posedness of the Spatially Homogeneous {{Landau}} Equation for Soft Potentials},
  author = {Guerin, H{\'e}l{\`e}ne and Fournier, Nicolas},
  year = {2008},
  month = jun,
  number = {arXiv:0806.3379},
  eprint = {0806.3379},
  primaryclass = {math},
  publisher = {arXiv},
  doi = {10.48550/arXiv.0806.3379},
  urldate = {2024-12-10},
  archiveprefix = {arXiv}
}

@misc{GuillenSilvestre2023,
  title = {The {{Landau}} Equation Does Not Blow Up},
  author = {Guillen, Nestor and Silvestre, Luis},
  year = {2023},
  month = nov,
  number = {arXiv:2311.09420},
  eprint = {2311.09420},
  primaryclass = {math},
  publisher = {arXiv},
  doi = {10.48550/arXiv.2311.09420},
  urldate = {2024-01-21},
  archiveprefix = {arXiv}
}

@article{HendersonSnelson2019,
  title = {C-Infinity {{Smoothing}} for {{Weak Solutions}} of the {{Inhomogeneous Landau Equation}}},
  author = {Henderson, Christopher and Snelson, Stanley},
  year = {2019},
  month = nov,
  publisher = {SPRINGER},
  issn = {0003-9527},
  doi = {10.1007/s00205-019-01465-7},
  urldate = {2025-02-26},
  copyright = {http://rightsstatements.org/vocab/InC/1.0/},
  langid = {english}
}

@misc{HendersonSnelsonTarfulea2019,
  title = {Local {{Solutions}} of the {{Landau Equation}} with {{Rough}}, {{Slowly Decaying Initial Data}}},
  author = {Henderson, Christopher and Snelson, Stanley and Tarfulea, Andrei},
  year = {2019},
  month = sep,
  number = {arXiv:1909.05914},
  eprint = {1909.05914},
  primaryclass = {math},
  publisher = {arXiv},
  doi = {10.48550/arXiv.1909.05914},
  urldate = {2025-01-27},
  archiveprefix = {arXiv}
}

@article{HendersonWang2024,
  title = {Kinetic {{Schauder}} Estimates with Time-Irregular Coefficients and Uniqueness for the {{Landau}} Equation},
  author = {Henderson, Christopher and Wang, Weinan},
  year = {Mon Apr 01 04:00:00 UTC 2024},
  journal = {Discrete and Continuous Dynamical Systems},
  volume = {44},
  number = {4},
  pages = {1026--1072},
  publisher = {{Discrete and Continuous Dynamical Systems}},
  issn = {1078-0947},
  doi = {10.3934/dcds.2023137},
  urldate = {2025-03-24},
  copyright = {http://creativecommons.org/licenses/by/3.0/},
  langid = {english}
}

@article{ImbertMouhot2021,
  title = {{The Schauder estimate in kinetic theory with application to a toy nonlinear model}},
  author = {Imbert, Cyril and Mouhot, Cl{\'e}ment},
  year = {2021},
  journal = {Annales Henri Lebesgue},
  volume = {4},
  pages = {369--405},
  issn = {2644-9463},
  doi = {10.5802/ahl.75},
  urldate = {2025-02-26},
  langid = {french}
}

@misc{Ji2024,
  title = {Dissipation Estimates of the {{Fisher}} Information for the {{Landau}} Equation},
  author = {Ji, Sehyun},
  year = {2024},
  month = oct,
  number = {arXiv:2410.09035},
  eprint = {2410.09035},
  publisher = {arXiv},
  urldate = {2024-11-12},
  archiveprefix = {arXiv}
}

@misc{Ji2024b,
  title = {Entropy Dissipation Estimates for the {{Landau}} Equation with {{Coulomb}} Potentials},
  author = {Ji, Sehyun},
  year = {2024},
  month = aug,
  number = {arXiv:2305.09841},
  eprint = {2305.09841},
  primaryclass = {math},
  publisher = {arXiv},
  doi = {10.48550/arXiv.2305.09841},
  urldate = {2025-02-11},
  archiveprefix = {arXiv}
}

@inproceedings{Krylov1996,
  title = {Lectures on {{Elliptic}} and {{Parabolic Equations}} in {{H{\"o}lder Spaces}}},
  author = {Krylov, N.},
  year = {1996},
  month = sep,
  series = {Graduate {{Studies}} in {{Mathematics}}},
  volume = {12},
  publisher = {American Mathematical Society},
  address = {Providence, Rhode Island},
  doi = {10.1090/gsm/012},
  urldate = {2025-02-26},
  isbn = {978-0-8218-0569-5 978-1-4704-2070-3},
  langid = {english}
}

@book{LifsicPitaevskij2008,
  title = {Physical Kinetics},
  editor = {Lif{\v s}ic, Evgenij M. and Pitaevskij, Lev P.},
  year = {2008},
  series = {Course of Theoretical Physics},
  number = {v. 10},
  publisher = {Elsevier},
  address = {Amsterdam Boston},
  isbn = {978-0-08-057049-5},
  langid = {english}
}

@article{NicolasFournierArnaudGuillin2017,
  title = {From a {{Kac-like}} Particle System to the {{Landau}} Equation for Hard Potentials and {{Maxwell}} Molecules},
  author = {{Nicolas Fournier} and {Arnaud Guillin}},
  year = {2017},
  journal = {Annales scientifiques de l'{\'E}cole normale sup{\'e}rieure},
  volume = {50},
  number = {1},
  pages = {157--199},
  issn = {0012-9593, 1873-2151},
  doi = {10.24033/asens.2318},
  urldate = {2024-09-03},
  langid = {english}
}

@misc{Silvestre2016,
  title = {Upper Bounds for Parabolic Equations and the {{Landau}} Equation},
  author = {Silvestre, Luis},
  year = {2016},
  month = aug,
  number = {arXiv:1511.03248},
  eprint = {1511.03248},
  primaryclass = {math},
  publisher = {arXiv},
  doi = {10.48550/arXiv.1511.03248},
  urldate = {2025-01-06},
  archiveprefix = {arXiv}
}

@article{Villani1998,
  title = {On a {{New Class}} of {{Weak Solutions}} to the {{Spatially Homogeneous Boltzmann}} and {{Landau Equations}}},
  author = {Villani, C{\'e}dric},
  year = {1998},
  month = sep,
  journal = {Archive for Rational Mechanics and Analysis},
  volume = {143},
  number = {3},
  pages = {273--307},
  issn = {1432-0673},
  doi = {10.1007/s002050050106},
  urldate = {2025-02-03},
  langid = {english}
}

@misc{Wiedemann2017,
  title = {Weak-Strong Uniqueness in Fluid Dynamics},
  author = {Wiedemann, Emil},
  year = {2017},
  month = may,
  number = {arXiv:1705.04220},
  eprint = {1705.04220},
  primaryclass = {math},
  publisher = {arXiv},
  doi = {10.48550/arXiv.1705.04220},
  urldate = {2025-02-13},
  archiveprefix = {arXiv}
}

\end{document}